\documentclass[10pt,reqno]{article}
\usepackage{amssymb, amsmath, amsthm, amsfonts, amscd}
\usepackage[english]{babel}
\usepackage{graphicx, epsfig}
\newtheorem{theorem}{Theorem}[section]
\newtheorem{cor}[theorem]{Corollary}
\newtheorem{lemma}[theorem]{Lemma}
\newtheorem{prop}[theorem]{Proposition}
\newtheorem{definition}{Definition}

\newcommand{\nm}{\noalign{\smallskip}}

\def\ep{\epsilon}

\newcommand{\Bx}{\mathbf{x}}
\newcommand{\By}{\mathbf{y}}
\newcommand{\Bk}{\mathbf{k}}

\newcommand{\RR}{\mathbb{R}}

\newcommand{\ZZ}{\mathbb{Z}}
\newcommand{\Scal}{\mathcal{S}}
\newcommand{\p}{\partial}

\newcommand{\pd}[2]{\frac {\p #1}{\p #2}}
\newcommand{\ds}{\displaystyle}
\newcommand{\eqnref}[1]{(\ref {#1})}
\newcommand{\pf}{\medskip \noindent {\sl Proof}. \ }
\renewcommand{\qed}{\hfill $\Box$ \medskip}
\newcommand{\beq}{\begin{equation}}
\newcommand{\eeq}{\end{equation}}

\numberwithin{equation}{section}
\numberwithin{figure}{section}

\begin{document}

\title{Enhancement of near-cloaking. Part II: the Helmholtz equation\thanks{\footnotesize
This work was supported  by ERC Advanced Grant Project
MULTIMOD--267184 and National Research Foundation through grants
No. 2010-0017532 and 2010-0004091.}}

\author{Habib Ammari\thanks{\footnotesize Department of Mathematics and Applications, Ecole Normale Sup\'erieure,
45 Rue d'Ulm, 75005 Paris, France (habib.ammari@ens.fr).} \and
Hyeonbae Kang\thanks{Department of Mathematics, Inha University,
Incheon 402-751, Korea (hbkang@inha.ac.kr, hdlee@inha.ac.kr).}
\and Hyundae Lee\footnotemark[3] \and Mikyoung
Lim\thanks{\footnotesize Department of Mathematical Sciences,
Korean Advanced Institute of Science and Technology, Daejeon
305-701, Korea (mklim@kaist.ac.kr). }}

\maketitle

\begin{abstract}
The aim of this paper is to extend the method of \cite{cloak_cond}
to scattering problems. We construct very effective near-cloaking
structures for the scattering problem at a fixed frequency. These
new structures are, before using the transformation optics,
layered structures and are designed so that their first scattering
coefficients vanish. Inside the cloaking region, any target has
near-zero scattering cross section for a band of frequencies. We
analytically show that our new construction significantly enhances
the cloaking effect for the Helmholtz equation.

\end{abstract}

\noindent {\footnotesize {\bf AMS subject classifications.} 35R30,
35B30}

\noindent {\footnotesize {\bf Key words.} cloaking, transformation
optics, Helmholtz equation, scattering cross section, scattering
coefficients}

\section{Introduction}
The cloaking problem in electromagnetic wave scattering is to make
a target invisible from far-field measurements. The difficulty is
to construct permeability and permittivity distributions of a
cloaking structure such that any target placed inside the
structure has zero scattering cross section. Since the pioneer
works \cite{glu}, \cite{leonhardt}, and \cite{pendry}, extensive
research has been done on cloaking in electromagnetic scattering.
We refer to \cite{GKLU}, \cite{sirev2}, and \cite{crossect} for
recent development on electromagnetic cloaking. One of the main
tools to obtain cloaking is to use a change of variables scheme
(also called transformation optics). The change of variables based
cloaking method uses a singular transformation to boost the
material property so that it makes a cloaking region look like a
point to outside measurements. However, this transformation
induces the singularity of material constants in the transversal
direction (also in the tangential direction in two dimensions),
which causes difficulty both in the theory and applications. To
overcome this weakness, so called `near cloaking' is naturally
considered, which is a regularization or an approximation of
singular cloaking. In \cite{kohn1}, instead of the singular
transformation, the authors use a regular one to push forward the
material constant in the conductivity equation, in which a small
ball is blown up to the cloaking region. In \cite{kohn2}, this
regularization point of view is adopted for the Helmholtz
equation. See also \cite{liu, nguyen}. It is worth mentioning that
there
  is yet another kind of cloaking in which  the cloaking region is outside the cloaking device, for instance,
  anomalous localized resonance \cite{MN_PRSA_06, MNMP_PRSA_05, inv}. See also \cite{alu}.

The purpose of this paper is to propose a new cancellation
technique in order to achieve enhanced invisibility from
scattering cross section.

Our approach extends to scattering problems the method first
provided in \cite{cloak_cond} to achieve near-cloaking for the
conductivity problem from Dirichlet-to-Neumann measurements. It is
based on the multi-coating which cancels the scattering
coefficients of the cloak. We first design a structure coated
around a perfect insulator to have vanishing scattering
coefficients of lower orders and show that the order of  the
scattering cross section of a small perfect insulator can be
reduced significantly. We then obtain near-cloaking structure by
pushing forward the multi-coated structure around a small object
via the standard blow-up transformation. Note that such structure
achieves near-cloaking for a band of frequencies. In order to
illustrate the viability of our method, we give numerical values
for the permittivity and permeability parameters and the thickness
of the layers of scattering coefficients vanishing structures.

This paper is organized as follows. In the next section we derive
the multi-polar expansion of the solution to the Helmholtz
equation, and define the scattering coefficients. In Section 3 we
characterize the scattering coefficients vanishing structures. In
Section 4 we show that the near-cloaking is enhanced  if a
scattering coefficients vanishing structure is used. In Section 5
we present some numerical examples of the scattering coefficients
vanishing structures.

Even though we consider only two dimensional Helmholtz equation in
this paper, the same argument can be applied to the equation in
three dimensions. The multi-coating technique developed in this
paper can be extended to the full Maxwell equations, which will be a subject of a forthcoming paper.

\section{Scattering coefficients}\label{subscat}

For $k >0$,   the fundamental solution to the Helmholtz operator
$(\Delta +k^2)$ in two dimensions satisfying
 $$
 (\Delta +k^2)\Gamma_k (\Bx) = \delta_0
 $$
with the outgoing radiation condition is given by
 $$
 \Gamma_k(\Bx) =- \frac{i}{4}H^{(1)}_0(k|\Bx|),
 $$
where $H_0^{(1)}$ is the Hankel function of the first kind of
order zero. For $|\Bx|>|\By|$, we recall Graf's addition formula
\cite{watson},
\begin{equation}\label{Graf}
H_0^{(1)}(k|\Bx-\By|) = \sum_{n \in \ZZ}
H_n^{(1)}(k|\Bx|)e^{in\theta_\Bx} J_n(k|\By|)e^{-in\theta_\By},
\end{equation}
where $\Bx=(|\Bx|,\theta_\Bx)$ and $\By=(|\By|,\theta_\By)$ in polar coordinates. Here $H_n^{(1)}$ is the
Hankel function of the first kind of order $n$ and $J_n$ is the Bessel function of order $n$.


We first define the scattering coefficients of an inclusion (with
two phase electromagnetic materials) and derive some important properties of them.

Let $D$ be a bounded domain in $\RR^2$ with Lipschitz boundary
$\partial D$, and let $(\epsilon_0, \mu_0)$ be the pair of
electromagnetic parameters (permittivity and permeability) of
$\RR^2\setminus\bar{D}$ and $(\epsilon_1, \mu_1)$ be that of $D$.
Then the permittivity and permeability distributions are given by
 \beq
 \ep = \epsilon_0\chi(\RR^2\setminus\bar{D})+\epsilon_1\chi(D) \quad\mbox{and}\quad \mu= \mu_0 \chi(\RR^2\setminus\bar{D})+ \mu_1 \chi(D).
 \eeq
Here and throughout this paper $\chi(D)$ is the characteristic function of $D$.

Given a frequency $\omega$, set $k= \omega \sqrt{\epsilon_1
\mu_1}$ and $k_0=\omega \sqrt{\epsilon_0 \mu_0}$. For a function
$U$ satisfying $(\Delta+k_0^2)U=0$ in $\RR^2$, we  consider the
scattered wave $u$, {\it i.e.}, the solution to the following
equation:
 \beq\label{HP1}
 \ \left \{
 \begin{array}{l}
 \ds\nabla \cdot \frac{1}{\mu} \nabla u
 +\omega^2 \ep u=0\quad \mbox{in } \RR^2, \\
 \nm \ds (u-U) \mbox{ satisfies the outgoing radiation condition}.
 \end{array}
 \right .
 \eeq

Let $\Scal_D^k [\varphi]$ be the single layer potential defined by the kernel $\Gamma_k$, {\it i.e.},
 $$
 \Scal_D^k [\varphi](\Bx)= \int_{\p D} \Gamma_k(\Bx-\By) \varphi(\By) d \sigma(\By).
 $$
Let $\Scal_D^{k_0}$ be the single layer potential associated with
the kernel $\Gamma_{k_0}$. Then, from, for example, \cite{AK04},
we know that the solution to \eqnref{HP1} can be represented using
the single layer potentials $\Scal_D^{k_0}$ and $\Scal_D^{k}$ as
follows
 \beq\label{helmScal}
 u(\Bx) = \begin{cases}
 U(\Bx)+\Scal_D^{k_0}[\psi](\Bx),\quad \Bx\in\RR^2\setminus\bar{D},\\
 \Scal_D^k [\varphi](\Bx),\quad \Bx\in D,
 \end{cases}
 \eeq
where the pair $(\varphi,\psi)\in L^2(\p D)\times L^2(\p D)$ is
the unique solution to \beq\label{phipsi} \ \left\{
\begin{array}{l}
\ds \Scal^k_D[\varphi] - \Scal_D^{k_0}[\psi] = U\\
\nm
\ds \left.\frac{1}{\mu}\pd{(\Scal_D^k[\varphi])}{\nu}\right|_-
-\left.\frac{1}{\mu_0}\pd{(\Scal_D^{k_0}[\psi])}{\nu}\right|_+=\frac{1}{\mu_0}\pd{U}{\nu}
\end{array}\right.
\quad\mbox{on }\p D. \eeq  Here, $+$ and $-$ in the subscripts
respectively indicate the limit from outside $D$ and inside $D$ to
$\p D$ along the normal direction and $\partial/\partial \nu$
denotes the normal derivative. It is proved in \cite{AK04} that
there exists a constant $C=C(k,k_0,D)$ such that
 \beq\label{primaryest}
 \|\varphi\|_{L^2(\p D)}+\|\psi\|_{L^2(\p D)}\leq C (\|U\|_{L^2(\p D)} +\|\nabla U\|_{L^2(\p
 D)}).
 \eeq
It is also proved in the same paper that there are constants
$\rho_0$ and $C=C(k,k_0,D)$ independent of $\rho$ as long as $\rho
\le \rho_0$ such that
 \beq\label{secondest}
 \|\varphi_\rho \|_{L^2(\p D)}+\|\psi_\rho \|_{L^2(\p D)}\leq C (\|U\|_{L^2(\p D)} +\|\nabla U\|_{L^2(\p
 D)}),
 \eeq
 where $(\varphi_\rho, \psi_\rho)$ is the solution of \eqnref{phipsi} with $k$ and $k_0$ respectively replaced by $\rho k$ and $\rho k_0$.

Note that the following asymptotic formula holds as $|\Bx| \to
\infty$, which can be seen from \eqnref{helmScal} and  Graf's
formula \eqnref{Graf}: \beq\label{helmfar}
u(\Bx)-U(\Bx)=-\frac{i}{4}\sum_{n\in\ZZ}H_n^{(1)}(k_0|\Bx|)e^{in\theta_\Bx}\int_{\p
D}J_n(k_0|\By|)e^{-in\theta_\By}\psi(\By)d\sigma(\By). \eeq Let
$(\varphi_m,\psi_m)$ be the solution to \eqnref{phipsi} with
$J_m(k_0|\Bx|)e^{im\theta_\Bx}$ in the place of $U(\Bx)$. We
define the {\it scattering coefficient} as follows.

\begin{definition} The scattering coefficients $W_{nm}$, $m,n \in \ZZ$, associated with the permittivity and
permeability distributions $\ep, \mu$ and the frequency $\omega$
(or $k, k_0, D$) are defined by
 \beq\label{stdef}
 W_{nm}= W_{nm}[\ep, \mu, \omega]:=\int_{\p D}J_n(k_0|\By|)e^{-in\theta_\By}\psi_{m}(\By)d\sigma(\By).
 \eeq
\end{definition}

We first obtain the following lemma for the size of $|W_{nm}|$.
\begin{lemma}\label{W_bound}
There is a constant $C$ depending on $(\ep, \mu, \omega)$ such
that \beq\label{wnmsize} |W_{nm}[\ep, \mu, \omega]|\leq
\frac{C^{|n|+|m|}}{|n|^{|n|} |m|^{|m|}} \quad \mbox{for all } n, m
\in \ZZ . \eeq Moreover, there exists $\rho_0$ such that, for all
$\rho \le \rho_0$,
 \beq\label{wnmsize2}
 |W_{nm}[\ep, \mu, \rho \omega]|\leq \frac{C^{|n|+|m|}}{|n|^{|n|} |m|^{|m|}}\rho^{|n|+|m|} \quad
 \mbox{for all }  n, m \in \ZZ,
 \eeq
where the constant $C$ depends on $(\ep, \mu, \omega)$ but is independent of $\rho$.
\end{lemma}

\pf Let $U(\Bx)= J_m(k_0|\Bx|)e^{im\theta_\Bx}$ and $(\varphi_m,\psi_m)$ be the solution to \eqnref{phipsi}. Since
 \beq\label{largebessel}
 J_m(t)\sim\frac{1}{\sqrt{2\pi |m|}}\Bigr(\frac{et}{2|m|}\Bigr)^{|m|}
 \eeq
as $m\rightarrow\infty$ (see \cite{AS}), we have
 $$
 \|U\|_{L^2(\p D)} +\|\nabla U\|_{L^2(\p D)} \le \frac{C^{|m|}}{|m|^{|m|}}
 $$
for some constant $C$. Thus it follows from \eqnref{primaryest} that
 \beq\label{psimest}
 \| \psi_m \|_{L^2(\p D)} \le \frac{C^{|m|}}{|m|^{|m|}}
 \eeq
for another constant $C$. So we get \eqnref{wnmsize} from \eqnref{stdef}.

Now let $(\varphi_m,\psi_m)$ be the solution to \eqnref{phipsi} with $k$ and $k_0$ replaced by $\rho k$ and $\rho k_0$.   One can see from \eqnref{secondest} that \eqnref{psimest} still holds for some $C$ independent of $\rho$ as long as $\rho \le \rho_0$ for some $\rho_0$.  Note that
 $$
 W_{nm}[\ep, \mu, \rho\omega] =\int_{\p D}J_n(\rho
 k_0|\By|)e^{-in\theta_\By}\psi_{m,\rho}(\By)d\sigma(\By),
 $$
where $(\varphi_{m,\rho}, \psi_{m,\rho})$ is the solution to
\eqnref{phipsi} with $k$ and $k_0$ respectively replaced by $\rho
k$ and $\rho k_0$ and $J_m(k_0\rho |\Bx|) e^{im\theta_\Bx}$ in the
place of $U(\Bx)$. So one can use \eqnref{largebessel} to obtain
\eqnref{wnmsize2}. This completes the proof. \qed

If $U$ is given as
 \beq\label{generalU}
 U(\Bx)=\sum_{m\in\ZZ}a_m(U) J_m(k_0|\Bx|)e^{im\theta_\Bx}
 \eeq
where $a_m(U)$ are constants, it follows from the principle of superposition that the solution $(\varphi, \psi)$ to \eqnref{phipsi} is given by
 $$
 \psi= \sum_{m\in \ZZ} a_m(U) \psi_m.
 $$
Then one can see from \eqnref{helmfar} that the solution $u$ to \eqnref{HP1} can be represented as
 \beq\label{Helmrep}
 u(\Bx)-U(\Bx)  = -\frac{i}{4}\sum_{n\in\ZZ}H_n^{(1)}(k_0|\Bx|)e^{in\theta_\Bx}\sum_{m\in\ZZ}W_{nm}a_{m}(U)
 \quad\mbox{as }|\Bx|\rightarrow\infty.
 \eeq
In particular, if $U$ is given by a plane wave $e^{i{\bf{k}}\cdot
\Bx}$ with  ${\bf{k}}\cdot{\bf{k}}=k_0^2$, then
 \beq\label{planewave}
 u(\Bx)-e^{i{\bf{k}}\cdot \Bx}=-\frac{i}{4}\sum_{n\in\ZZ}H_n^{(1)}(k_0|\Bx|)e^{in\theta_\Bx}
 \sum_{m\in\ZZ}W_{nm}e^{im(\frac{\pi}{2}-\theta_{\bf{k}})} \quad\mbox{as }|\Bx|\rightarrow\infty,
 \eeq
where
${\bf{k}}=k_0(\cos\theta_{\bf{k}},\sin\theta_{\bf{k}})$ and
$\bf{x}=(|\bf{x}|, \theta_{\Bx})$. In fact, since
$$
e^{ik_0r\sin\theta}=\sum_{m\in\ZZ} J_m(k_0 r)e^{im\theta},
$$
we have
\begin{equation}\label{planewaveexp}
e^{i{\bf{k}}\cdot \Bx} =\sum_{m\in\ZZ}
e^{im(\frac{\pi}{2}-\theta_{\bf{k}})}J_m(k_0|\Bx|)e^{im\theta_\Bx},
\end{equation}
and
 \beq\label{planepsi}
 \psi= \sum_{m\in \ZZ} e^{im(\frac{\pi}{2}-\theta_{\bf{k}})} \psi_m.
 \eeq
Thus \eqnref{planewave} holds. It is worth emphasizing that the
expansion formula \eqnref{Helmrep} or \eqnref{planewave}
determines uniquely the scattering coefficients $W_{nm}$, for $n,
m\in\ZZ$.

We now show that the scattering coefficients are basically the
Fourier coefficients of the far-field pattern (the scattering
amplitude) which is $2\pi$-periodic function in two dimensions.
The far-field pattern $A_\infty [\epsilon, \mu, \omega]$, when the
incident field is given by $e^{i{\bf{k}}\cdot \Bx}$, is defined to
be
 \beq\label{umU}
 u(\Bx)-e^{i{\bf{k}}\cdot \Bx}  = -i e^{-\frac{\pi i}{4}} \frac{ e^{ik_0|\Bx|} }{\sqrt{8\pi k_0|\Bx|}} A_\infty [\epsilon, \mu, \omega](\theta_{\Bk},\theta_\Bx) + o(|\Bx|^{-\frac{1}{2}})
  \quad\mbox{as }|\Bx|\rightarrow\infty.
 \eeq

Recall that
 \beq\label{hbehavior}
 H_0^{(1)}(t)\sim \sqrt{\frac{2}{\pi t}}e^{i(t-\frac{\pi}{4})} \quad\mbox{as } t \rightarrow\infty,
 \eeq
where $\sim$ indicates that the difference between the right-hand
and left-hand side is $O(t^{-1})$.  If $|\Bx|$ is large while
$|\By|$ is bounded, then we have
$$|\Bx - \By|= |\Bx|-|\By|\cos(\theta_\Bx-\theta_\By)+O(\frac{1}{|\Bx|}),$$
and hence
 $$
 H_0^{(1)}(k_0 |\Bx-\By|)  \sim e^{-\frac{\pi i}{4}} \sqrt{\frac{2}{\pi k_0 |\Bx|}} e^{ik_0 (|\Bx| - |\By|\cos(\theta_\Bx - \theta_\By))} \quad \mbox{as } |\Bx| \to \infty.
 $$
Thus, from \eqnref{helmScal}, we get
 \beq
 u(\Bx)-e^{i{\bf{k}}\cdot \Bx} \sim -ie^{-\frac{\pi i}{4}} \frac{ e^{ik_0|\Bx|} }{\sqrt{8\pi k_0|\Bx|}}
 \int_{\p D}e^{-ik_0|\By|\cos(\theta_\Bx - \theta_\By)}\psi(\By)\ d\sigma(\By) \quad\mbox{as }|\Bx|\rightarrow\infty,
 \eeq
and infer that the far-field pattern is given by
 \beq\label{scatterform}
 A_\infty(\theta_\Bk, \theta_\Bx) = \int_{\p D}e^{-ik_0|\By|\cos(\theta_\Bx - \theta_\By)}\psi(\By)\ d\sigma(\By),
 \eeq
where $\psi$ is given by \eqnref{planepsi}.

Let
 $$
 A_\infty(\theta_\Bk,\theta_\Bx)=\sum_{n\in\ZZ}b_n (\theta_\Bk) e^{in\theta_\Bx}
 $$
be the Fourier series of $A_\infty(\theta_\Bk,\cdot)$. From
\eqnref{scatterform} it follows that
 \begin{align*}
 b_n(\theta_\Bk) &= \frac{1}{2\pi} \int_{0}^{2\pi} \int_{\p D} e^{-ik_0|\By|\cos(\theta_\Bx - \theta_\By)} \psi(\By)\, d\sigma(\By) \, e^{-in\theta_\Bx} \, d\theta_\Bx \\
 &= \frac{1}{2\pi} \int_{\p D} \int_{0}^{2\pi} e^{-ik_0|\By|\cos(\theta_\Bx - \theta_\By)} e^{-in\theta_\Bx} \, d\theta_\Bx \, \psi(\By) \, d\sigma(\theta_\By).
 \end{align*}
Since
 $$
 \frac{1}{2\pi} \int_{0}^{2\pi} e^{-ik_0|\By|\cos(\theta_\Bx - \theta_\By)} e^{-in\theta_\Bx} \, d\theta_\Bx = J_n (k_0 |\By|) e^{-in (\theta_\By + \frac{\pi}{2})},
 $$
we deduce that
 $$
 b_n(\theta_\Bk) = \int_{\p D} J_n (k_0 |\By|) e^{-in (\theta_\By + \frac{\pi}{2})} \psi(\By) \, d\sigma(\theta_\By).
 $$
Using \eqnref{planepsi} we now arrive at the following theorem.
\begin{theorem} \label{2DFT} Let $\theta$ and $\theta'$ be  respectively the incident
and scattered direction. Then we have
 \beq \label{Ainfty}
 A_\infty [\epsilon, \mu, \omega](\theta, \theta') =  \sum_{n, m\in\ZZ} i^{(m-n)} e^{in\theta'} W_{nm} [\epsilon, \mu, \omega] e^{-i m\theta}.
 \eeq
\end{theorem}
We emphasize that the series in \eqnref{Ainfty} converges uniformly thanks to \eqnref{wnmsize}. We also note that
if $U$ is given by \eqnref{generalU} then the scattering amplitude, which we denote by $A_\infty [\epsilon, \mu, \omega](U, \theta')$, is given by
 \beq\label{farexp}
 A_\infty [\epsilon, \mu, \omega](U, \theta')= \sum_{n\in\ZZ}
 i^{-n}  e^{in\theta'}\sum_{m\in\ZZ}W_{nm}a_{m}(U).
 \eeq
The conversion of the far-field to the near field is achieved via
formula \eqnref{planewave}.

The scattering cross section $S[\epsilon, \mu, \omega]$ is defined by
 \beq
 S[\epsilon, \mu, \omega](\theta'):= \int_0^{2\pi} \bigg|A_\infty
 [\epsilon, \mu, \omega](\theta, \theta') \bigg|^2\, d\theta.
 \eeq
See \cite{born,taylor}. As an immediate consequence of Theorem \ref{2DFT} we obtain the following corollary.
\begin{cor} \label{2DFT2}
We have
 \beq\label{scsform}
 S[\epsilon, \mu, \omega](\theta') = 2\pi \sum_{m\in\ZZ} \bigg|\sum_{n\in\ZZ} i^{-n} W_{nm} [\epsilon, \mu, \omega]
 e^{in\theta'} \bigg|^2.
 \eeq
\end{cor}

It is worth mentioning that the optical theorem \cite{born,taylor} leads to a
natural constraint on $W_{nm}$.  In fact, we have
 \beq\label{constraint}
 \Im m \, A_\infty[\epsilon, \mu, \omega](\theta', \theta') = - \sqrt{\frac{\omega}{8\pi}}
 S[\epsilon, \mu, \omega](\theta'),\quad \forall \; \theta' \in [0, 2\pi],
 \eeq
or equivalently,
 \beq \label{constraint2}
 \Im m \, \sum_{n, m\in\ZZ} i^{m-n} e^{i(n-m)\theta'} W_{nm} [\epsilon,
 \mu, \omega]  = - \sqrt{\frac{\pi \omega}{2}} \sum_{m\in\ZZ}
 \bigg|\sum_{n\in\ZZ} i^{-n} W_{nm} [\epsilon, \mu, \omega]
 e^{in\theta'} \bigg|^2,
 \quad \forall \; \theta' \in [0, 2\pi].
 \eeq

In the next section, we compute the scattering coefficients of
multiply coated inclusions and provide structures whose scattering
coefficients vanish. Such structures will be used in Section
\ref{sect4} to enhance near cloaking. Any target placed inside
such structures will have nearly vanishing scattering cross
section $S$, uniformly in  the direction $\theta^\prime$.

\section{S-vanishing structures}

The purpose of this section is to construct multiply layered
structures whose scattering coefficients vanish. We call such
structures {\it S-vanishing structures}. We design a multi-coating
around an insulated inclusion $D$, for which the scattering
coefficients vanish. The computations of the scattering
coefficients of multi-layered structures (with multiple phase
electromagnetic materials) follow in exactly the same way as in
Section \ref{subscat}. The system of two equations
\eqnref{phipsi} should be replaced by a system of $2\times$ the
number of phase interfaces ($-1$ if the core is perfectly
insulating).

For positive numbers $r_1,\ldots, r_{L+1}$ with
$2=r_1>r_2>\cdots>r_{L+1}=1$, let
 $$
 A_j:=\{x: r_{j+1} \le |x| < r_j\}, \quad j=1,\dots,L, \quad A_0:=\RR^2 \setminus
 \overline{A_1}, \quad A_{L+1} (=D):=\{x:  |x| < 1\}.
 $$
Let $(\mu_j, \epsilon_j)$ be the pair of permeability and permittivity of $A_j$ for
$j=0,1,\ldots, L+1$. Set $\mu_0=1$ and $\epsilon_0=1$. Let
 \beq\label{STstructure2}
 \mu = \sum_{j=0}^{L+1} \mu_j \chi(A_j)
 \quad \mbox{and} \quad \epsilon = \sum_{j=0}^{L+1} \epsilon_j
 \chi(A_j).
 \eeq

In this case the scattering coefficient
$W_{nm}=W_{nm}[\mu,\epsilon,\omega]$ can be defined using
\eqnref{Helmrep}. In fact, if $u$ is the solution to
 \beq\label{muepU}
 \nabla\cdot\frac{1}{\mu}\nabla u+\omega^2\epsilon u=0\qquad\mbox{in }\RR^2
 \eeq
with the outgoing radiation condition on $u-U$ where $U$ is given
by \eqnref{generalU}, then $u-U$ admits the asymptotic expansion
\eqnref{Helmrep} with $k_0= \omega \sqrt{\epsilon_0 \mu_0}$.

Exactly like the conductivity case \cite{cloak_cond} one can show
using symmetry that \beq W_{nm}=0\quad\mbox{if }m\neq n. \eeq Let
us define $W_n$ by \begin{equation} \label{defwn} W_n:=W_{nn}.
\end{equation}
Our purpose is to design,  given $N$ and $\omega$, $\mu$ and
$\epsilon$ so that $W_{n}[\mu,\epsilon,\omega]=0$ for $|n| \le N$.
We call such a structure $(\mu, \epsilon)$ an \emph{S-vanishing
structure of order $N$ at frequency $\omega$}. Since
$H^{(1)}_{-n}=(-1)^n H^{(1)}_{n}$ and $J_{-n}=(-1)^n J_{n}$, we
have
 \beq
 W_{-n,-n}= {W}_{nn},
 \eeq
and hence it suffices to consider $W_{nn}$ only for $n \ge 0$.

Note that \eqnref{constraint2} leads to
 \beq \label{constraint3}
 \Im m \, \sum_{n \in\ZZ} W_{n} [\epsilon, \mu, \omega]  = -
 \sqrt{\frac{\pi \omega}{2}}  \sum_{n \in\ZZ} \bigg| W_{n} [\epsilon, \mu, \omega] \bigg|^2.
 \eeq

Let $k_j :=\omega\sqrt{\mu_j\epsilon_j}$ for $j=0,1, \ldots, L$. We assume that  $\mu_{L+1}=\infty$, which amounts to that the solution satisfies the zero Neumann condition on $|\Bx|=r_{L+1} (=1)$.
To compute $W_n$ for $n \ge 0$, we look for solutions $u_n$ to
\eqnref{muepU} of the form
\begin{equation}\label{Helsol}
u_n(\Bx) = a_j^{(n)} J_n(k_j r)
e^{in\theta}+b_j^{(n)}H_n^{(1)}(k_j r)e^{in\theta},\quad \Bx \in
A_j, \quad j=0,\dots,L,
\end{equation}
with $a_0^{(n)}=1$. Note that
 \beq
 W_n =4ib_0^{(n)}.
 \eeq
The solution $u_n$ satisfies the transmission conditions
 $$
 u_n |_{+}=u_n|_{-} \quad \mbox{and}\quad \frac{1}{\mu_{j-1}} \frac{\p u_n}{\p \nu} \Big|_{+} = \frac{1}{\mu_{j}} \frac{\p u_n}{\p \nu} \Big|_{-} \quad\mbox{on } |\Bx|=r_j
 $$
for $j=1, \ldots, L$, which reads
\begin{align}
&
\begin{bmatrix}
\ds   J_n(k_j r_j) & \ds H_n^{(1)}(k_j r_j)\\
\ds \sqrt{\frac{\epsilon_j}{\mu_j}} J'_n(k_j r_j) &
\ds\sqrt{\frac{\epsilon_j}{\mu_j}} \left(H^{(1)}_n\right)'(k_j
r_j)
\end{bmatrix}
\begin{bmatrix}
\ds a_j^{(n)}\\
\ds b_j^{(n)}
\end{bmatrix} \nonumber \\
&=\begin{bmatrix}
\ds   J_n(k_{j-1} r_j) & \ds H_n^{(1)}(k_{j-1} r_j)\\
\ds\sqrt{\frac{\epsilon_{j-1}}{\mu_{j-1}}}  J'_n(k_{j-1} r_j) &
\ds\sqrt{\frac{\epsilon_{j-1}}{\mu_{j-1}}}
\left(H^{(1)}_n\right)'(k_{j-1} r_j)
\end{bmatrix}
\begin{bmatrix}
\ds a_{j-1}^{(n)}\\
\ds b_{j-1}^{(n)}
\end{bmatrix}. \label{tranmatrix}
\end{align}
The Neumann condition $\frac{\p u_n}{\p \nu}|_{+}=0$ on $|\Bx|=r_{L+1}$ amounts to
\beq\label{neumanmatrix}
\begin{bmatrix} 0&0\\ J'_n( k_L) & \left(H^{(1)}_n\right)'( k_L)
\end{bmatrix}
\begin{bmatrix}
\ds a_{L}^{(n)}\\
\ds b_{L}^{(n)}
\end{bmatrix} = \begin{bmatrix}
\ds0\\
\ds0
\end{bmatrix} .
\eeq
Combining \eqnref{tranmatrix} and \eqnref{neumanmatrix}, we obtain
\beq\label{anplusone}
\begin{bmatrix}
\ds 0\\ \ds 0
\end{bmatrix}
= P^{(n)}[\epsilon, \mu, \omega]
\begin{bmatrix}
\ds a_0^{(n)}\\ \ds b_0^{(n)}
\end{bmatrix},
\eeq where
\begin{align*}&P^{(n)}[\epsilon, \mu, \omega]:= \begin{bmatrix}
\ds 0& \ds 0\\
\ds p^{(n)}_{21}& \ds p^{(n)}_{22}
\end{bmatrix}=(-\frac{\pi}{2} i \omega)^{L}\left(
\prod_{j=1}^{L} \mu_j r_j\right) \begin{bmatrix} 0&0\\ J'_n( k_L)
 & \left(H^{(1)}_n\right)'( k_L) \end{bmatrix}\\
&\times \prod_{j=1}^{L} \begin{bmatrix}
\ds\sqrt{\frac{\epsilon_j}{\mu_j}} \left(H^{(1)}_n\right)'( k_j
r_j)
& \ds - H_n^{(1)}( k_j r_j)\\
\nm \ds - \sqrt{\frac{\epsilon_j}{\mu_j}} J'_n( k_j r_j) & \ds
J_n( k_j r_j)
\end{bmatrix}
 \begin{bmatrix}
\ds   J_n( k_{j-1} r_j) & \ds H_n^{(1)}( k_{j-1} r_j)\\
\nm \ds \sqrt{\frac{\epsilon_{j-1}}{\mu_{j-1}}} J'_n( k_{j-1} r_j)
&
\ds\sqrt{\frac{\epsilon_{j-1}}{\mu_{j-1}}}\left(H^{(1)}_n\right)'(
k_{j-1} r_j)
\end{bmatrix}.\end{align*}

In order to have a structure whose scattering coefficients $W_n$ vanishes up to
the $N$-order, we need to have $b_0^{(n)}=0$ (when $a_0^{(n)}=1$) for $n=0,\ldots, N$, which amounts to
 \beq\label{ptwoone}
 p^{(n)}_{21} =0 \quad\mbox{for } n=0,\ldots, N
 \eeq
because of \eqnref{anplusone}. We emphasize that $p^{(n)}_{22}\ne0$. In fact, if $p^{(n)}_{22}=0$, then \eqnref{anplusone} can be fulfilled with $a_0^{(n)}=0$ and $b_0^{(n)}=1$. It means that there exists $(\mu, \ep)$ on $\RR^2 \setminus D$ such that the following problem has a solution:
 \beq\label{muepU2}
 \begin{cases}
 \ds \nabla\cdot\frac{1}{\mu}\nabla u+\omega^2\epsilon u=0\qquad\mbox{in }\RR^2 \setminus \overline{D},\\
 \nm
 \ds \pd{u}{\nu} \Big|_{+} =0 \quad\mbox{on } \p D, \\
 \nm
 u(\Bx)= H^{(1)}_n(k_0r)e^{in\theta}  \quad\mbox{for } |\Bx|=r > 2,
 \end{cases}
 \eeq
which is not possible.

We note that \eqnref{ptwoone} is a set of conditions on $(\mu_j, \ep_j)$ and $r_j$ for $j=1, \ldots, L$. In fact, $p^{(n)}_{21}$ is a nonlinear algebraic function of $\mu_j$, $\ep_j$ and $r_j$, $j=1, \ldots, L$. We are not able to show existence of $(\mu_j, \ep_j)$ and $r_j$, $j=1, \ldots, L$, satisfying \eqnref{ptwoone} even if it is quite important to do so. But the solutions (at fixed frequency) can be computed numerically in the same way as in the conductivity case \cite{cloak_cond}.

We now consider the S-vanishing structure for all (low) frequencies. Let $\omega$ be fixed and we look for a structure $(\mu, \ep)$ such that
 \beq
 W_n[\mu,\epsilon,\rho\omega]=0 \quad \mbox{for all } |n| \le N \mbox{ and } \rho \le \rho_0
 \eeq
for some $\rho_0$. Such a structure may not exist. So instead we
look for a structure such that
 \beq \label{wnhd}
 W_n[\mu,\epsilon,\rho\omega]=o(\rho^{2N}) \quad \mbox{for all } |n| \le N \mbox{ and } \rho \to 0.
 \eeq
We call such a structure an \emph{S-vanishing structure of order $N$ at low frequencies}.

To investigate the behavior of $W_n[\mu,\epsilon,\rho\omega]$ as
$\rho \to 0$, let us recall the behavior of Bessel functions for
small arguments. As $t \to 0$, it is known that
\begin{align*}
J_n(t)&=\frac{t^n}{2^n}\left(\frac{1}{\Gamma(n+1)}-\frac{\frac{1}{4}t^2}{\Gamma(n+2)}
+\frac{(\frac{1}{4}t^2)^2}{2!\Gamma(n+3)}-\frac{(\frac{1}{4}t^2)^3}{3!\Gamma(n+4)}+\cdots\right),\\
Y_n(t)&=-\frac{(\frac{1}{2}t)^{-n}}{\pi}\sum_{l=0}^{n-1}\frac{(n-l-1)!}{l!}(\frac{1}{4}t^2)^l
+\frac{2}{\pi}\ln(\frac{1}{2}t)J_n(t)\\
& \quad -\frac{(\frac{1}{2}t)^n}{\pi}\sum_{l=0}^\infty(\psi(l+1)+\psi(n+l+1))\frac{(-\frac{1}{4}t^2)^l}{l!(n+l)!},
\end{align*}
where $\psi(1)= -\gamma$ and $\psi(n)= -\gamma + \sum_{l=1}^{n-1}
1/l$ for $n\geq 2$ with $\gamma$ being the Euler constant.
In particular, if $n=0$, we have
\begin{align}
J_0(t)&=1-\frac{1}{4}t^2+\frac{1}{64}t^4+O(t^6), \\
Y_0(t)&=\frac{2}{\pi}\ln
t+\frac{2}{\pi}(\gamma-\ln2)-\frac{1}{2\pi}t^2\ln
t+\Bigr(\frac{1}{2\pi}-\frac{1}{2\pi}(\gamma-\ln2)\Bigr)t^2+O(t^4\ln
t).
\end{align}
Plugging these formulas into \eqnref{anplusone}, we have
\begin{align}
& P^{(0)}[\epsilon, \mu, \rho \omega]
 =(-\frac{\pi}{2}i  \rho \omega)^{L} \left( \prod_{j=1}^{L} \mu_j r_j \right)
 \begin{bmatrix} 0&0 \\
 \ds -\frac{k_L}{2} \rho +O(\rho^3) &
 \ds \frac{2i}{\pi k_L} \rho^{-1}+O(\rho \ln \rho)\end{bmatrix} \nonumber\\
& \nonumber \qquad \times \prod_{j=1}^{L}
\begin{bmatrix}
\ds \frac{2i}{\pi \omega\mu_j r_j}  \rho^{-1} + O(\rho\ln \rho) & \ds
\frac{4}{\pi^2}\left(\frac{1}{\omega \mu_{j-1} r_j} -
\frac{1}{\omega \mu_{j} r_j}\right)\frac{\ln \rho}{\rho}+O(\rho^{-1}) \\
\nm
\ds\frac{r_j}{2}\omega\epsilon_j\left(1-\frac{\epsilon_{j-1}}{
\epsilon_j}\right)\rho
+ O(\rho^3) & \ds \frac{2i}{\pi \omega\mu_{j-1} r_j}  \rho^{-1} + O(\rho\ln
\rho)
\end{bmatrix}\\
&=\ds \rho^{-1}\begin{bmatrix} 0&0\\  O(\rho^2)&\ds \frac{2i}{\pi
k_L}\prod_{j=1}^L \frac{\mu_j}{\mu_{j-1}}+O(\rho)
 \end{bmatrix}, \label{P0}
\end{align}
and, for $n\geq 1$,
\begin{align}
& P^{(n)}[\epsilon, \mu, \rho\omega] = (-i \frac{\pi}{2} \rho
\omega)^{L}\left( \prod_{j=1}^{L} \mu_j r_j \right)
 \begin{bmatrix} 0&0\\
 \ds \frac{nk_L^{n-1}}{2^n\Gamma(n+1)}\rho^{n-1} + O(\rho^{n})
 & \ds \frac{i2^n\Gamma(n+1)}{\pi k_L^{n+1}}  \rho^{-n-1}
 + O(\rho^{-n})\end{bmatrix} \nonumber\\
& \qquad\quad  \times \prod_{j=1}^{L}
\begin{bmatrix}
\ds \sqrt{\frac{\epsilon_j}{\mu_j}}\frac{i2^n\Gamma(n+1)}{\pi(k_j
r_j)^{n+1}}  \rho^{-n-1} + O(\rho^{-n})
& \ds \frac{i2^n\Gamma(n)}{\pi(k_j r_j)^n}\rho^{-n} + O(\rho^{-n+1}) \\
\nm \ds -\sqrt{\frac{\epsilon_j}{\mu_j}} \frac{n(k_j
r_j)^{n-1}}{2^n\Gamma(n+1)}\rho^{n-1} + O(\rho^{n}) & \ds \frac{(k_j
r_j)^n} {2^n\Gamma(n+1)} \rho^n + O(\rho^{n+1})
\end{bmatrix} \nonumber\\
& \nonumber \qquad\quad \times \begin{bmatrix} \ds   \frac{(k_{j-1}
r_j)^n}{2^n\Gamma(n+1)}\rho^n + O(\rho^{n+1})
& \ds -\frac{i2^n\Gamma(n)}{\pi(k_{j-1} r_j)^n}\rho^{-n} + O(\rho^{-n+1}) \\
\nm \ds \sqrt{\frac{\epsilon_{j-1}}{\mu_{j-1}}} \frac{n(k_{j-1}
r_j)^{n-1}}{2^n\Gamma(n+1)}\rho^{n-1} + O(\rho^{n}) & \ds
\sqrt{\frac{\epsilon_{j-1}}{\mu_{j-1}}}\frac{i2^n\Gamma(n+1)}{\pi(k_{j-1}
r_j)^{n+1}}  \rho^{-n-1} + O(\rho^{-n})
\end{bmatrix} \\
\nm & \nonumber \qquad = \frac{1}{2^L}
\begin{bmatrix} 0&0\\
\ds \frac{nk_L^{n-1}}{2^n\Gamma(n+1)}\rho^{n-1} + O(\rho^{n})
& \ds \frac{i2^n\Gamma(n+1)}{\pi k_L^{n+1}}  \rho^{-n-1} + O(\rho^{-n})\end{bmatrix}\\
&\qquad\quad \times\prod_{j=1}^{L}
\begin{bmatrix}
\ds  a_j (b_j+1)  + o(1)
& \ds  c_j ( b_j-1) \rho^{-2n} + o(\rho^{-2n}) \\
\nm \ds \frac{b_j-1}{c_j} \rho^{2n} + o(\rho^{2n}) & \ds
\frac{b_j+1}{a_j} + o(1)
\end{bmatrix}, \label{P1_first}
\end{align}
where
 $$
 a_j :=  \left( \frac{k_{j-1}}{k_j} \right)^n, \quad b_j := \frac{\mu_{j}}{\mu_{j-1}}, \quad c_j := \frac{i 2^{2n} \Gamma(n)\Gamma(n+1)}{\pi(k_{j-1} k_j r_j^2)^n},
 $$
with $k_j=\omega \sqrt{\epsilon_j \mu_j}$.

From the above calculations of the leading order terms of
$P^{(n)}[\epsilon, \mu, \rho\omega]$ and the expansion formula of
$J_n(t)$ and $Y_n(t)$, we see that $p_{21}^{(n)}$ and
$p_{22}^{(n)}$ admit the following expansions:
 \beq \label{pop21}
 p_{21}^{(n)}(\mu, \epsilon, t) = t^{n-1}
 \left(f_0^{(n)} (\mu,\epsilon)+\sum_{l=1}^{(N-n)}\sum_{j=0}^{L+1}
 f_{l,j}^{(n)} (\mu,\epsilon) t^{2l}(\ln t)^j+o(t^{2N-2n})\right)
 \eeq
 and
 \beq \label{pop22}
 p_{22}^{(n)}(\mu, \epsilon, t) =  t^{-n-1} \left(g_0^{(n)} (\mu,\epsilon)
 +\sum_{l=1}^{(N-n)}\sum_{j=0}^{L+1} g_{l,j}^{(n)} (\mu,\epsilon) t^{2l}(\ln t)^j+o(t^{2N-2n})\right) \eeq
for $t= \rho\omega$ and some functions $f_0^{(n)}, g_0^{(n)}, f_{l,j}^{(n)},$ and
$g_{l,j}^{(n)}$ independent of $t$.

\begin{lemma} For any pair of $(\mu, \epsilon)$, we have
 \beq
 g_0^{(n)} (\mu,\epsilon)\ne 0.
 \eeq
\end{lemma}
\pf For $n=0$, it follows from \eqnref{P0} that
$$  g_0^{(0)} (\mu,\epsilon)=\frac{2i}{\pi
\sqrt{ \epsilon_L\mu_L}}\prod_{j=1}^L \frac{\mu_j}{\mu_{j-1}}\ne 0.$$

Suppose $n>0$. Assume that there exists a pair of $(\mu,
\epsilon)$ such that $g^{(n)}_0(\mu,\epsilon)= 0.$ Then the
solution given by \eqnref{Helsol} with $a_0^{(n)}=0$ and
$b_0^{(n)}=1$ satisfies \beq\label{muepU4}
 \begin{cases}
 \ds \nabla\cdot\frac{1}{\mu}\nabla u+\rho^2\omega^2\epsilon u=0\qquad\mbox{in }\RR^2 \setminus \overline{D},\\
 \nm
 \ds \pd{u}{\nu} \Big|_{+} =o(\rho^{-n}) \quad\mbox{on } \p D, \\
 \nm
 u(\Bx)= H^{(1)}_n(\rho k_0r)e^{in\theta} \quad\mbox{for } |\Bx|=r > 2.
 \end{cases}
 \eeq
Let $v(\Bx):= \lim_{\rho\rightarrow 0} \rho^{n} u(\Bx)$. Then $v$ satisfies
\beq\label{muepV4}
 \begin{cases}
 \ds \nabla\cdot\frac{1}{\mu}\nabla v=0\qquad\mbox{in }\RR^2 \setminus \overline{D},\\
 \nm
 \ds \pd{v}{\nu} \Big|_{+} =0 \quad\mbox{on } \p D, \\
 \nm
 v(\Bx)=\ds -\frac{i2^n\Gamma(n)}{\pi(\epsilon_0 \mu_0)^{n/2}}r^{-n}e^{in\theta} \quad\mbox{for } |\Bx|=r > 2,
 \end{cases}
 \eeq
which is impossible. Thus $g_0^{(n)} (\mu,\epsilon)\ne 0,$ as
desired and the proof is complete. \qed

Equations \eqnref{pop21} and \eqnref{pop22} together with the above lemma
give us the following proposition.
\begin{prop} \label{prop3.1} For $n\ge 1$, let $W_n$ be defined by \eqnref{defwn}. We have
 \beq\label{homogen}
 W_n[\mu,\epsilon,t] (\theta,\theta') = t^{2n} \left({W}^{0}_n[\mu,\epsilon] (\theta,\theta') +
 \sum_{l=1}^{(N-n)}\sum_{j=0}^{M_{n,l}} {W}^{l,j}_n [\mu,\epsilon] (\theta,\theta') t^{2l}(\ln t)^j\right) + o(t^{2N})
 \eeq
uniformly in $(\theta,\theta')$, where $t=\rho\omega$, $M_{n,l}:(L+1) (N-n)$ ($L$ being the number of layers), and the
coefficients ${W}^{0}_n[\mu,\epsilon]$ and
${W}^{l,j}_n[\mu,\epsilon]$ are independent of $t$.
\end{prop}

To construct an S-vanishing structure of order $N$ at low
frequencies, we need to have a pair $(\mu,\epsilon)$ of the form
\eqnref{STstructure2} satisfying
 \beq\label{Wcond}
 {W}^{0}_n[\mu,\epsilon]=0,
 \mbox{ and } {W}^{l,j}_n [\mu,\epsilon] =0 \quad\mbox{for }
 0\leq n\leq N,\ 1\leq l\leq (N-n),\ 1\leq j\leq M_{n,l}.
 \eeq
We construct numerically such structures for small $N$ in the last section.

\section{Enhancement of near cloaking} \label{sect4}

In this section we show that the S-vanishing structures (after a
transformation optics) enhance the near cloaking.

Let $(\mu, \ep)$ be a S-vanishing structure of order $N$ at low frequencies, \emph{i.e.}, \eqnref{Wcond} holds, and it is of the form \eqnref{STstructure2}. It follows from \eqnref{wnmsize2}, Theorem \ref{2DFT}, and Proposition \ref{prop3.1} that
 \beq\label{without}
 A_\infty [\mu,\ep, \rho \omega ](\theta,\theta')=o(\rho^{2N})
 \eeq
uniformly in $(\theta,\theta')$ if $\rho \le \rho_0$ for some $\rho_0$.

Let
 \beq
 \Psi_{\frac{1}{\rho}}(\Bx)=\frac{1}{\rho}\Bx, \quad \Bx\in\RR^2.
 \eeq
Then we have
 \beq\label{scaleid}
 A_\infty\Bigr[\mu\circ\Psi_\frac{1}{\rho},\epsilon\circ\Psi_\frac{1}{\rho}, \omega\Bigr]
 =A_\infty [ \mu, \epsilon, \rho\omega ] .
 \eeq
To see this, let $u$ be the solution to \beq
 \ \left \{
\begin{array}{l}
\ds\nabla\cdot\frac{1}{(\mu\circ\Psi_{\frac{1}{\rho}})(\Bx)}\nabla u(\Bx)+\omega^2(\epsilon
\circ\Psi_\frac{1}{\rho})(\Bx) u(\Bx) = 0\quad \mbox{in } \RR^2 \setminus \overline{B_\rho},\\
\nm \ds
\pd{u}{\nu}=0 \quad\mbox{on } \partial B_\rho,\\
\nm \ds (u-U) \mbox{ satisfies the outgoing radiation condition},
\end{array}
\right .
 \eeq
where $U(\Bx) = e^{ik_0(\cos\theta,\sin\theta)\cdot \Bx}$. Here $B_\rho$ is the disk of radius $\rho$ centered at $0$. Define for $\By=\frac{1}{\rho} \Bx$
 $$
 \tilde{u}(\By):=\Bigr(u\circ\Psi_\frac{1}{\rho}^{-1}\Bigr)(\By) =\Bigr(u\circ\Psi_\rho\Bigr)(\By) \quad \mbox{and} \quad \tilde{U}(\By) =\Bigr(U\circ\Psi_\rho\Bigr)(\By).
 $$
Then, we have
\beq
 \ \left \{
\begin{array}{l}
\ds\nabla\cdot\frac{1}{\mu(\By)}\nabla_\By \tilde u(\By)+\rho^2\omega^2\epsilon(\By)\tilde u(\By) = 0\quad \mbox{in } \RR^2, \\
\nm
\ds \pd{\tilde u}{\nu}=0 \quad\mbox{on } \partial B_1,\\
\nm \ds (\tilde{u}-\tilde U) \mbox{ satisfies the outgoing
radiation condition}.
\end{array}
\right .
 \eeq
From the definition of the far-field pattern $A_\infty$, we get
 $$
 ( u- U)(\Bx)  \sim -ie^{-\frac{\pi i}{4}}\frac{ e^{ik_0 |\Bx|} }{
 \sqrt{8\pi k_0 |\Bx|}}A_\infty\Bigr[\mu\circ\Psi_\frac{1}{\rho},\epsilon\circ\Psi_\frac{1}{\rho},
 \omega\Bigr](\theta,\theta')
 \quad\mbox{as }|\Bx|\rightarrow\infty,
 $$
and
$$(\tilde u-\tilde U)(\By)  \sim -ie^{-\frac{\pi i}{4}}\frac{
e^{i\rho k_0 |\By|} }{\sqrt{8\pi \rho k_0 |\By|}}A_\infty[\mu,\epsilon,\rho\omega](\theta,\theta')
\quad\mbox{as }|\By|\rightarrow\infty,
$$
where $\Bx = |\Bx|(\cos\theta',\sin\theta')$. So, we have \eqnref{scaleid}. It then follows from \eqnref{without} that
 \beq\label{afterscale1}
 A_\infty\Bigr[\mu\circ\Psi_\frac{1}{\rho},\epsilon\circ\Psi_\frac{1}{\rho},
 \omega\Bigr](\theta,\theta')=o(\rho^{2N}).
 \eeq
We also obtain from \eqnref{scsform}
 \beq\label{afterscale2}
 S\Bigr[\mu\circ\Psi_\frac{1}{\rho},\epsilon\circ\Psi_\frac{1}{\rho},
 \omega\Bigr](\theta')=o(\rho^{4N}).
 \eeq
It is worth emphasizing that $(\mu\circ\Psi_\frac{1}{\rho},\epsilon\circ\Psi_\frac{1}{\rho})$ is a multi-coated structure of radius $2\rho$.

We now apply a transformation to the structure
$(\mu\circ\Psi_\frac{1}{\rho},\epsilon\circ\Psi_\frac{1}{\rho})$
to blow up the small disk of radius $\rho$. Let us recall the
following well-known lemma (see, for instance, \cite{GKLU}).
\begin{lemma}\label{scaling}
Let $F$ be a diffeomorphism of $\RR^2$ onto $\RR^2$ such that
$F(\Bx)$ is identity for $|\Bx|$ large enough. If $v$ is a
solution to
$$
\nabla \cdot A \nabla v + \omega^2q v=0 \quad \mbox{in } \RR^2,
$$
subject to the outgoing radiation condition, then $w$ defined by
$w(\By) = v(F^{-1}(\By))$ satisfies \beq \nabla \cdot (F_* A)
\nabla w + \omega^2 (F_* q) w =0 \quad\mbox{in } \RR^2, \eeq
together with the outgoing radiation condition, where \beq (F_*
A)(\By) = \frac{DF(\Bx)A(\Bx)DF^T(\Bx)}{\det DF(\Bx)} \quad
\mbox{and} \quad (F_* q)(\By) = \frac{q(\Bx)}{\det DF(\Bx)} \eeq
with $\Bx=F^{-1}(\By)$ and $T$ being the transpose.
\end{lemma}

For a small number $\rho$, let
$F_{\rho}$ be the diffeomorphism defined by
\begin{equation}
F_{\rho}(\Bx):= \begin{cases}
\ds \Bx \quad&\mbox{for } |\Bx|\geq2,\\
\nm \ds\Bigr(\frac{3-4\rho}{2(1-\rho)}+\frac{1}{4(1-\rho)}|\Bx|\Bigr)\frac{\Bx}{|\Bx|} \quad&\mbox{for }2\rho\leq|\Bx|\leq 2,\\
\nm \ds\Bigr(\frac{1}{2}+\frac{1}{2\rho}|\Bx|\Bigr)\frac{\Bx}{|\Bx|} \quad&\mbox{for }\rho\leq|\Bx|\leq 2\rho,\\
\nm \ds\frac{\Bx}{\rho} &\mbox{for }|\Bx|\leq \rho.
\end{cases}
\end{equation}
We then get from \eqnref{afterscale1}, \eqnref{afterscale2}, and Lemma \ref{scaling} the main result of this paper.
\begin{theorem} \label{thmmain}
If $(\mu, \ep)$ is a S-vanishing structure of order $N$ at low frequencies, then there exists $\rho_0$ such that
 \beq\label{aftertrans}
 A_\infty\Bigr[(F_{\rho})_*(\mu\circ\Psi_\frac{1}{\rho}),(F_{\rho})_*(\epsilon\circ\Psi_\frac{1}{\rho}),
 \omega\Bigr](\theta,\theta')=o(\rho^{2N}),
 \eeq
and
 \beq\label{aftertrans2}
 S\Bigr[(F_{\rho})_*(\mu\circ\Psi_\frac{1}{\rho}),(F_{\rho})_*(\epsilon\circ\Psi_\frac{1}{\rho}),
 \omega\Bigr](\theta')=o(\rho^{4N}),
 \eeq
for all $\rho \le \rho_0$, uniformly in $\theta$ and $\theta'$.
\end{theorem}
We first note the cloaking enhancement is achieved for all the
frequencies smaller than $\omega$. This is because \eqnref{wnhd}
holds if we replace $\omega$ by $\omega^\prime \leq \omega$. Then
it is worth comparing \eqnref{aftertrans} with \eqnref{without}.
In \eqnref{without}, $(\mu, \ep)$ is a multiply layered structure
between radius 1 and 2 in which each layer is filled with an
isotropic material, and enhanced near cloaking is achieved for low
frequencies $\rho \omega$ with $\rho \le \rho_0$. On the other
hand, in \eqnref{aftertrans} the frequency $\omega$ does not have
to be small. In fact, \eqnref{aftertrans} says that for any
frequency $\omega$ there is a radius $\rho$ which yields the
enhanced near cloaking up to $o(\rho^{2N})$. We emphasize that
$(F_{\rho})_*(\mu\circ\Psi_\frac{1}{\rho})$ and
$(F_{\rho})_*(\epsilon\circ\Psi_\frac{1}{\rho})$ are anisotropic
permittivity and permeability distributions. It is not clear
whether we can achieve enhanced near cloaking at high frequencies
by using isotropic layers as done in \eqnref{without}.

\section{Numerical examples}

In this section we provide numerical examples of S-vanishing
structures of order $N$ at low frequencies, \emph{i.e.},
structures $(\mu,\epsilon)$ of the form \eqnref{STstructure2}
satisfying \eqnref{Wcond}. To do so, we use the gradient descent
method to minimize over $(\mu,\epsilon)$ the quantity
 \beq\label{minfunc}
 |{W}^{0}_n[\mu,\epsilon] |^2 + \sum_{l=1}^{(N-n)}\sum_{j=0}^{M_{n,l}} | {W}^{l,j}_n [\mu,\epsilon]
 |^2,
 \eeq
where ${W}^{0}_n[\mu,\epsilon]$ and ${W}^{l,j}_n [\mu,\epsilon]$
are the coefficients of $W_n[\mu, \ep, t]$ in \eqnref{homogen}.
These coefficients are expressed in terms of Bessel functions and
their derivatives. We use the symbolic toolbox of MATLAB in order
to extract these coefficients.

It is quite challenging numerically to minimize the quantity in
\eqnref{minfunc} for large $N$. In numerical examples we take
$N=2$. In this case one can show through tedious computations that
the nonzero leading coefficients of $W_n[\mu, \epsilon, t]$ are as
follows:
\begin{itemize}
\item $[t^2,t^4, t^4\log t]$ for $n=0$, \emph{i.e.}, $W_0^{1,0}$,
$W_0^{2,0}$, $W_0^{2,1}$; \item $[t^2,t^4, t^4\log t]$ for $n=1$,
\emph{i.e.}, $W_1^{0}$, $W_1^{1,0}$, $W_1^{1,1}$; \item $[t^4]$
for $n=2$, \emph{i.e.}, $W_2^{0}$.
\end{itemize}

\begin{figure}[h!]
\begin{center}
\epsfig{figure=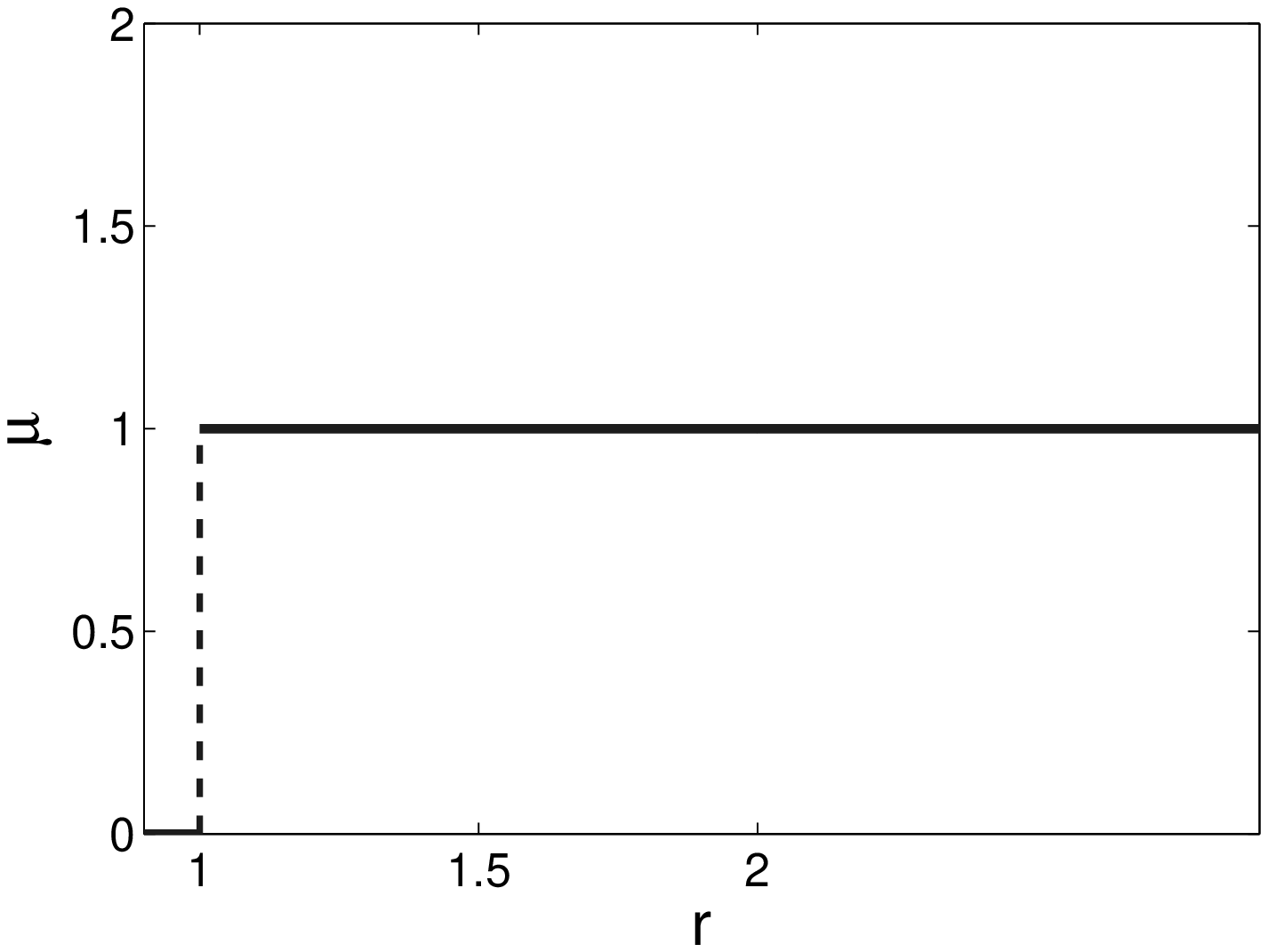, height=3.5cm}\hskip .2cm
\epsfig{figure=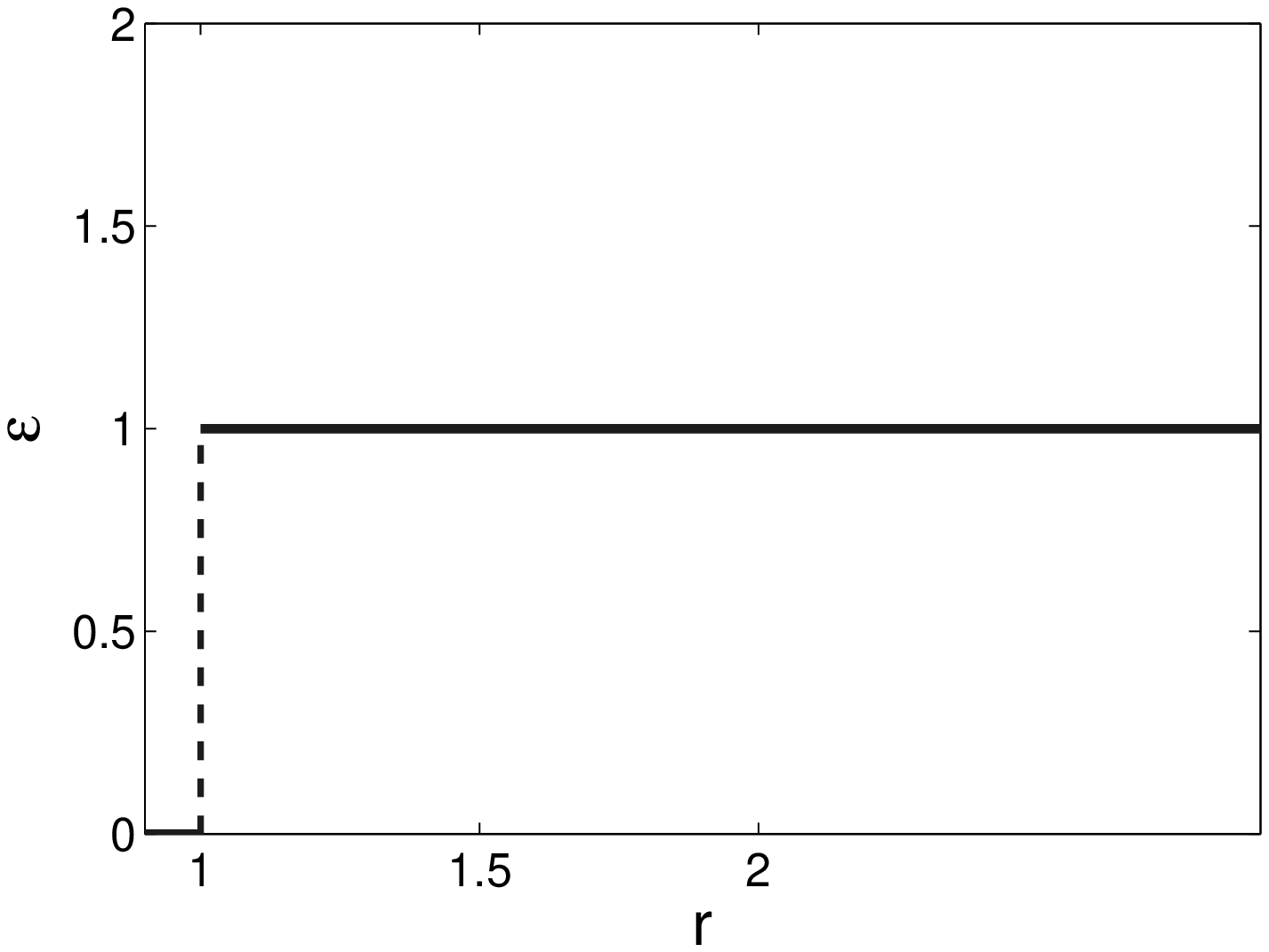,height=3.5cm}\hskip .2cm
\epsfig{figure=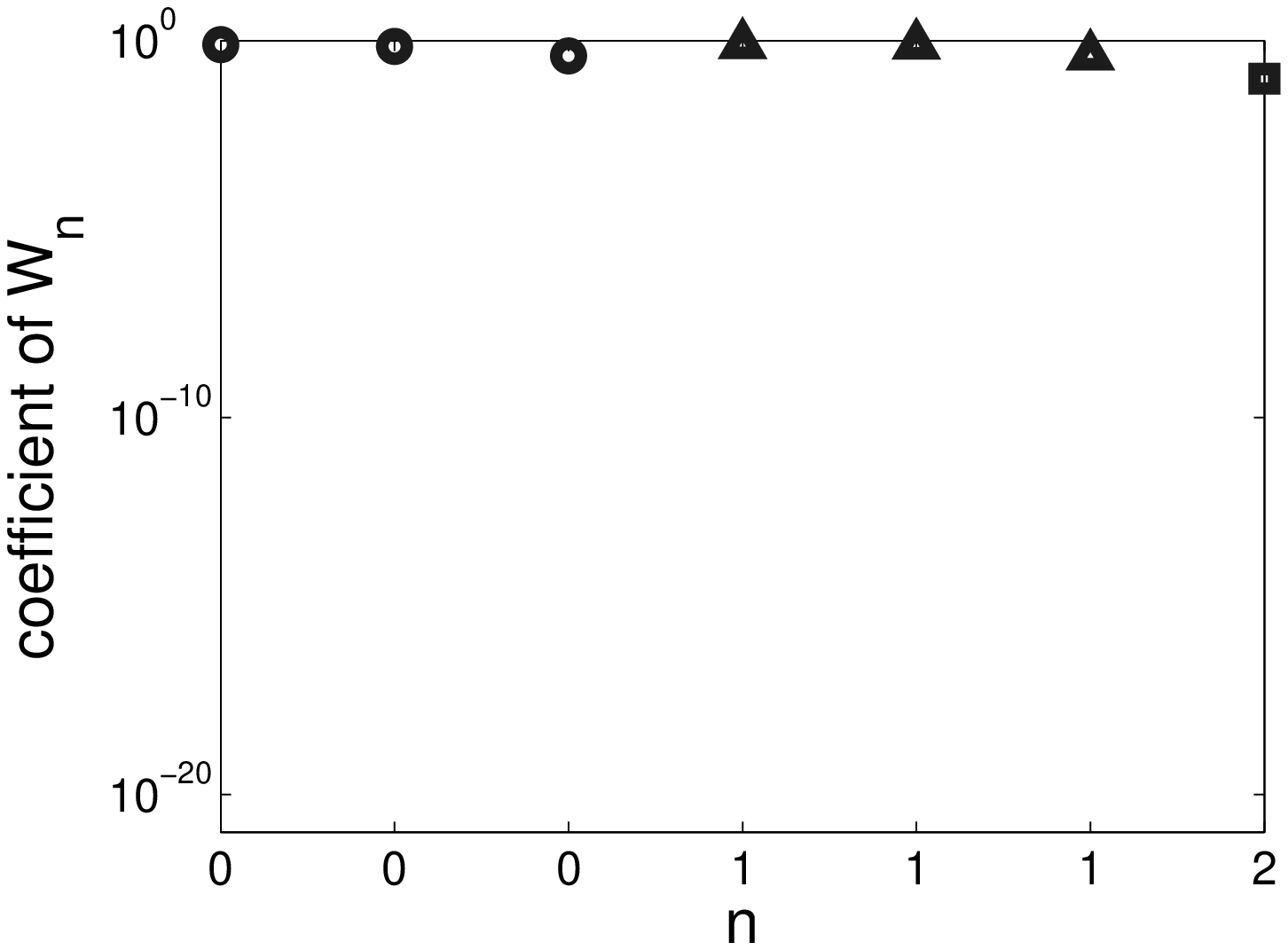,  height=3.5cm}\\
\epsfig{figure=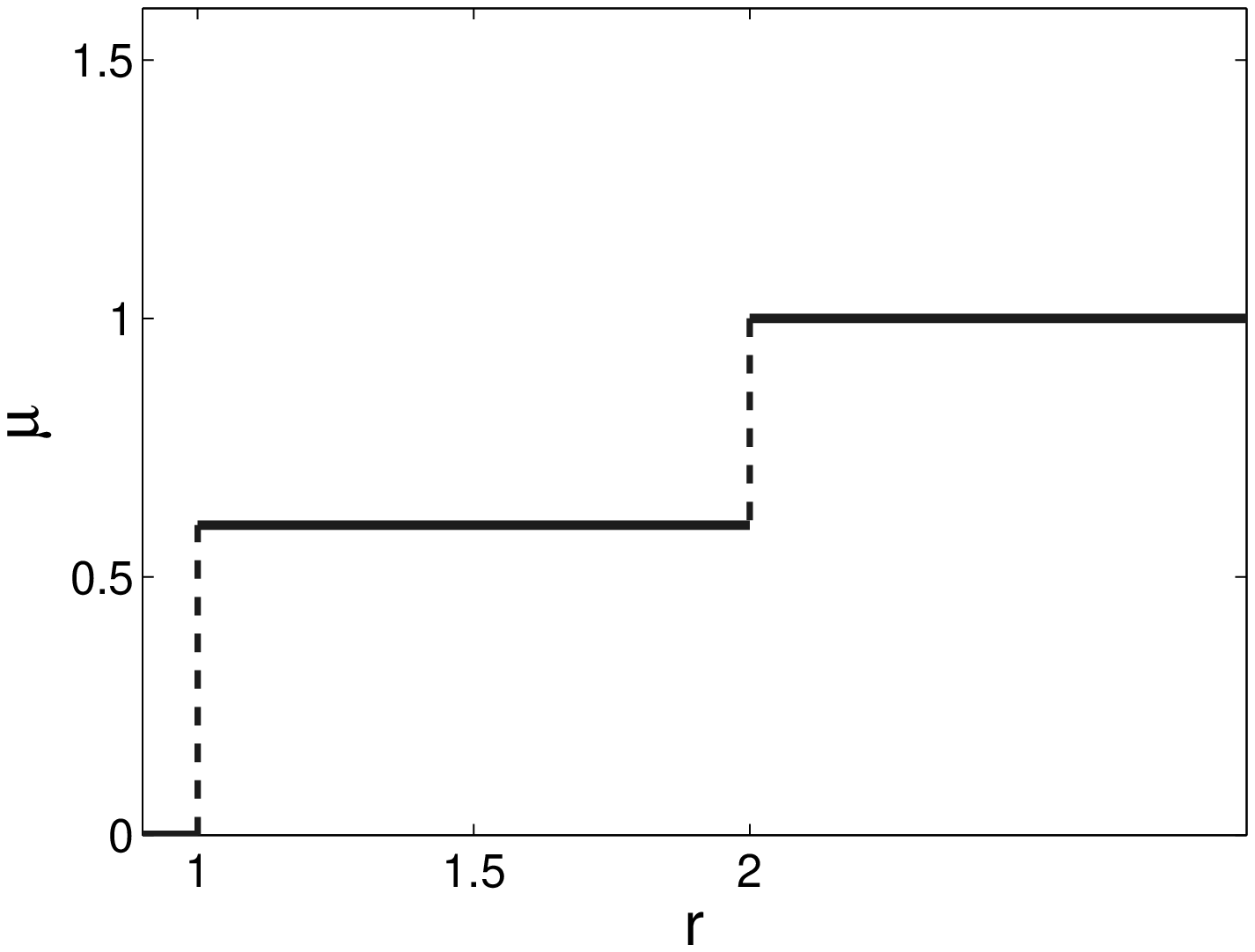, height=3.5cm}\hskip .2cm
\epsfig{figure=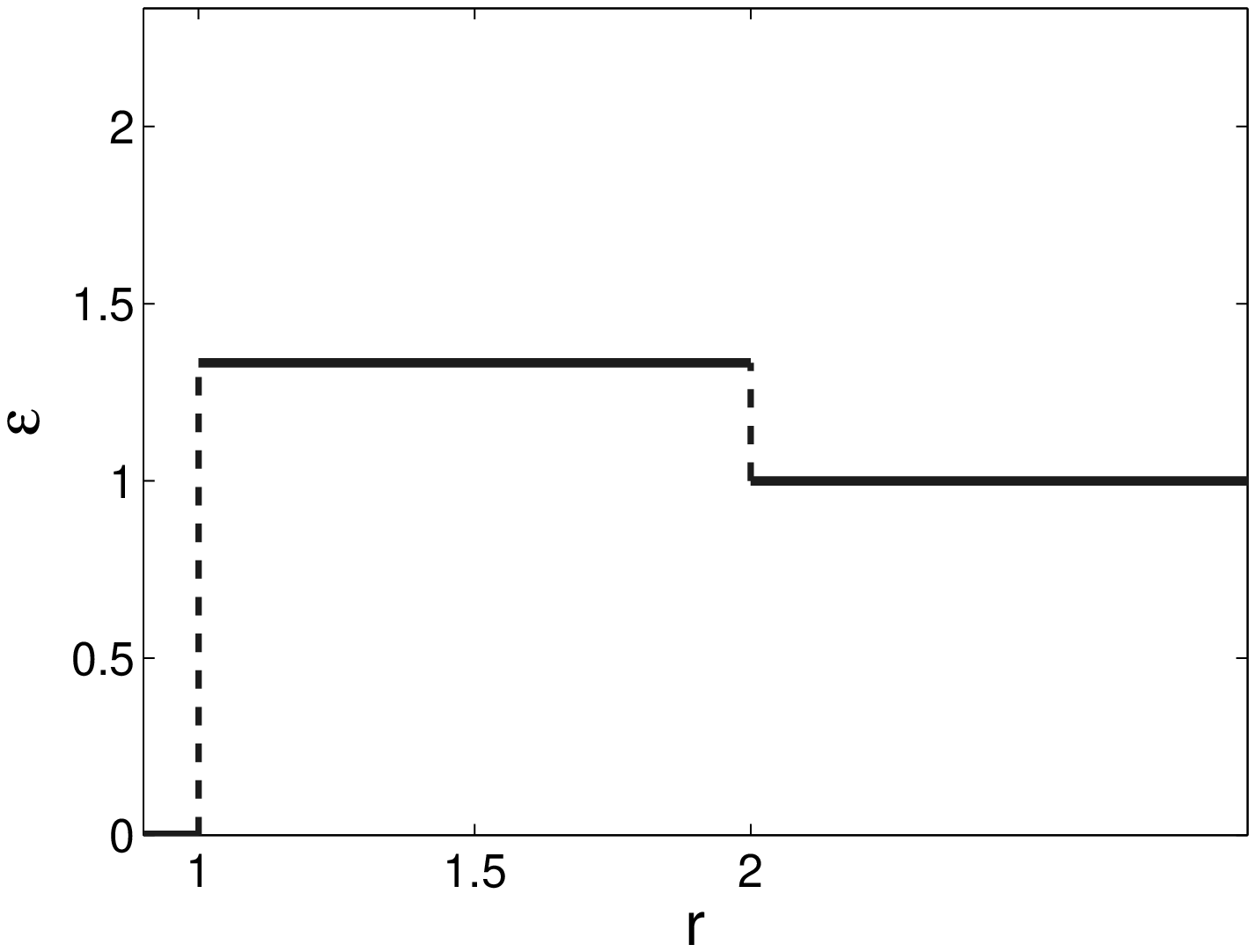,height=3.5cm}\hskip .2cm
\epsfig{figure=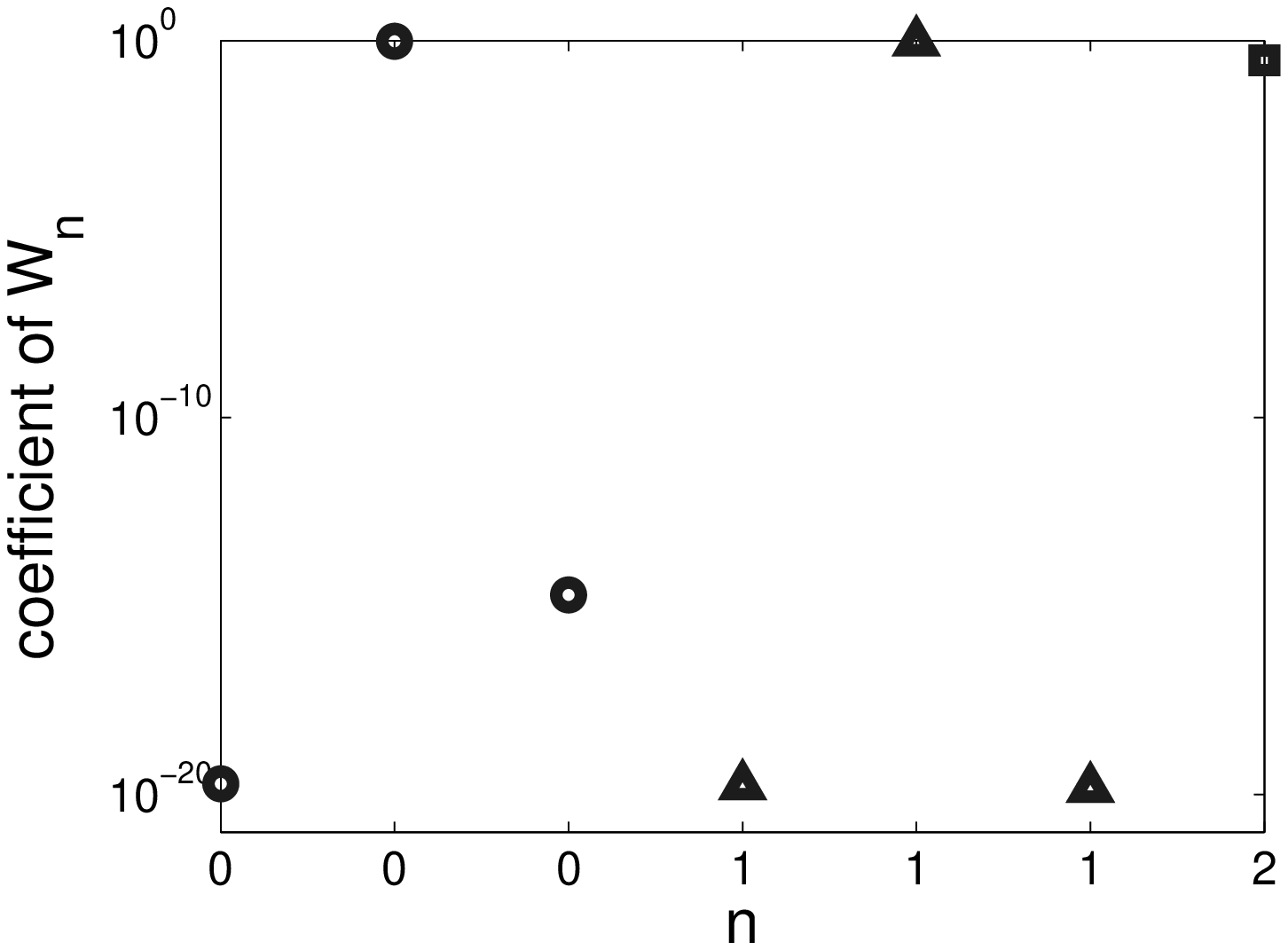,  height=3.5cm}\\
\epsfig{figure=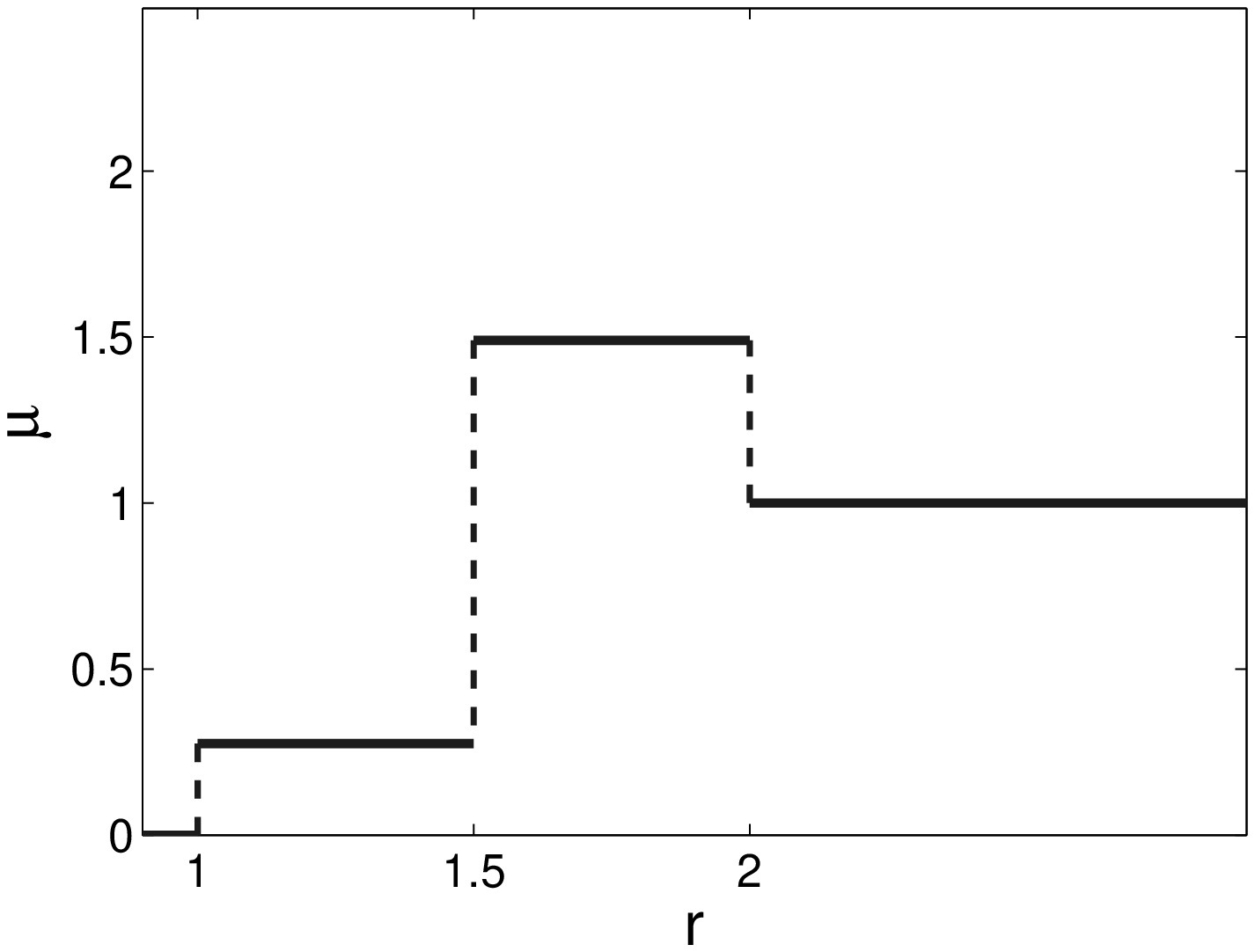, height=3.5cm}\hskip .2cm
\epsfig{figure=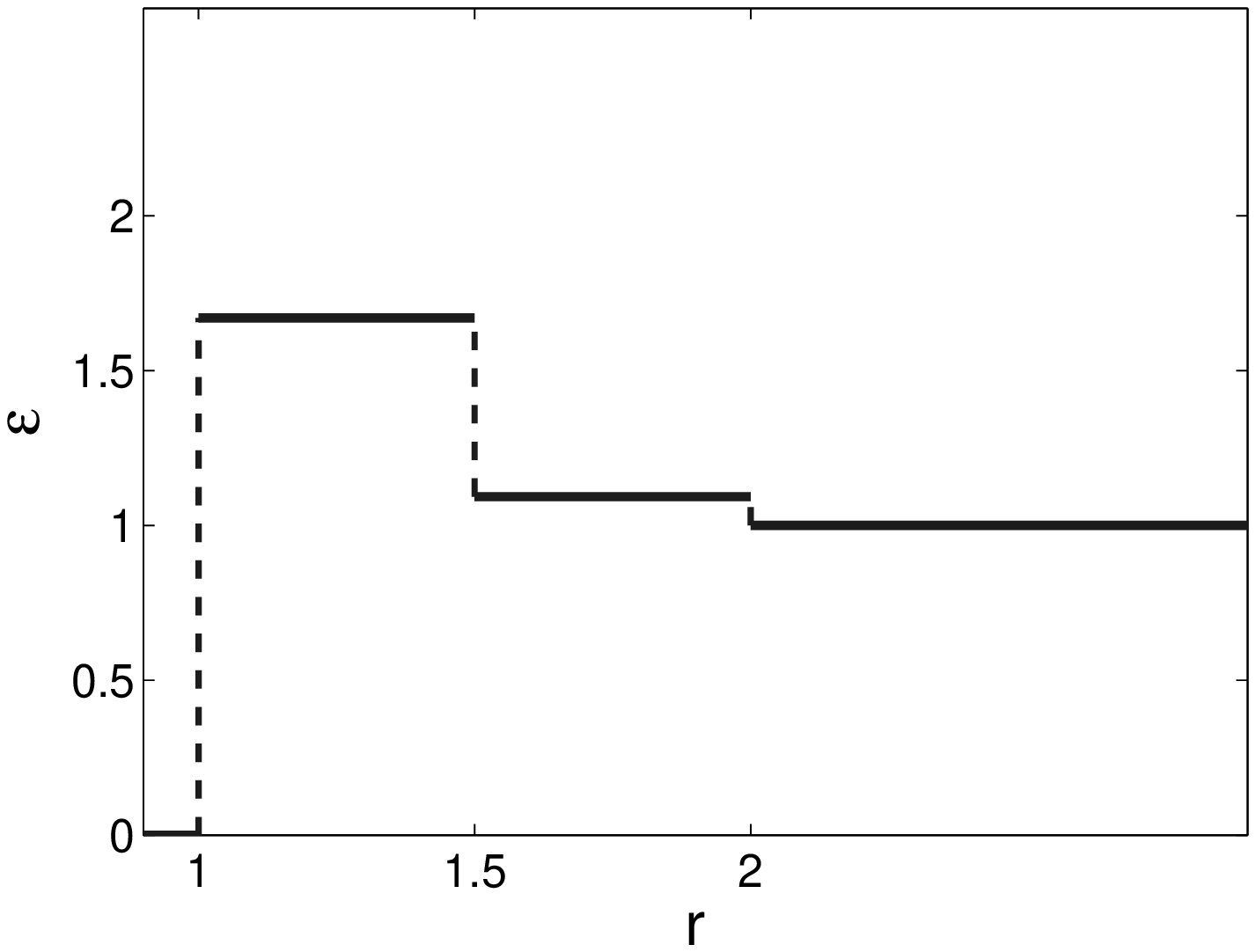,height=3.5cm}\hskip .2cm
\epsfig{figure=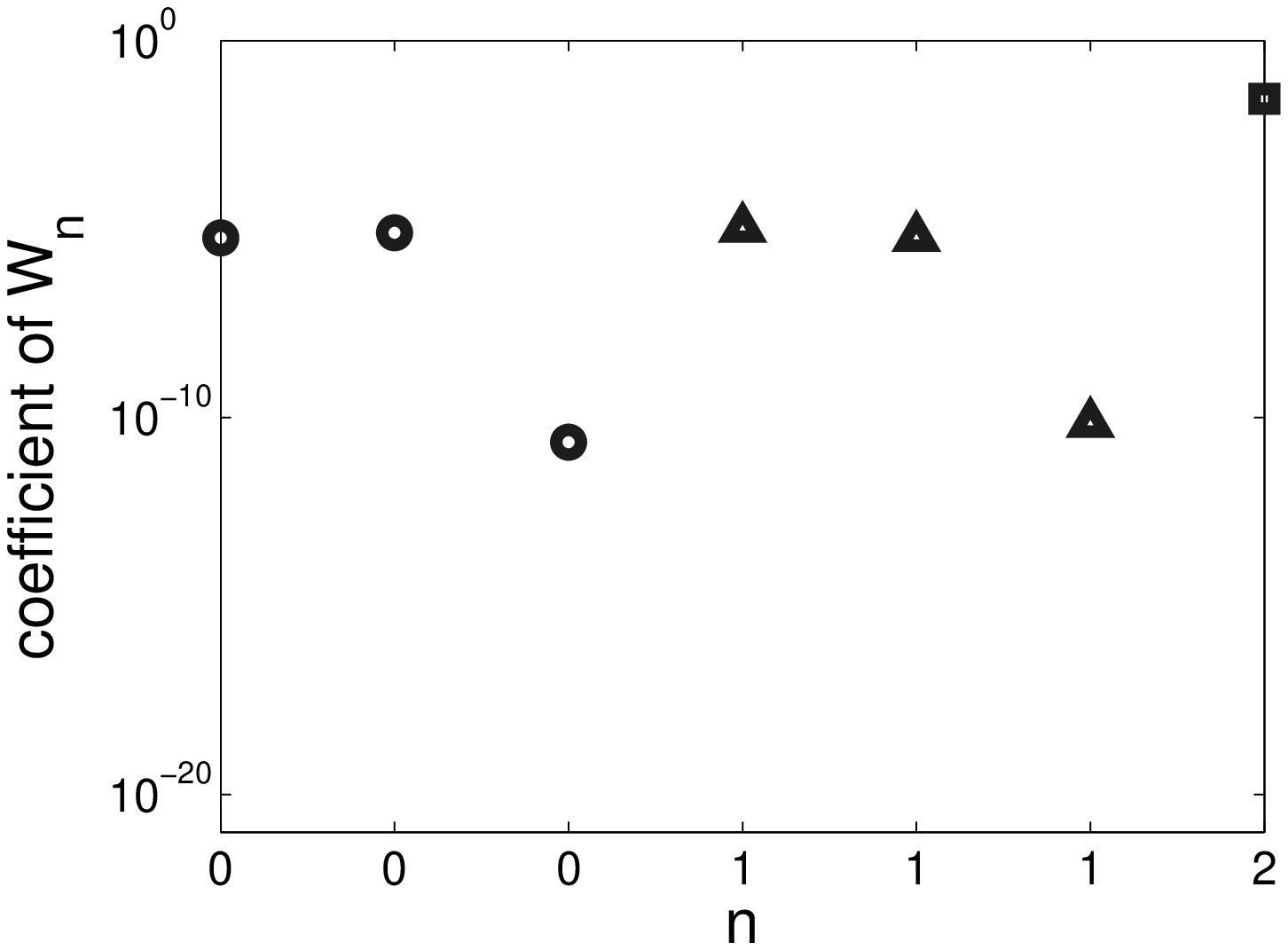,  height=3.5cm}\\
\end{center}
\caption{Graphs of the first and second columns show profiles of
permeability $\mu$ and the permittivity $\epsilon$.
The third column shows the coefficients of $[t^2,t^4, t^4\log t]$ in the expansion of
$W_0$ (represented by $(0,0,0)$) and $W_1$ (represented by
$(1,1,1)$), and the coefficient of $[t^4]$ in $W_2$ (represented by 2).
}\label{general}
\end{figure}

Results of computations are given in Figure \ref{general} when we
use 0, 1, 2 layers $(L=0,1,2)$. The computed material parameters
are $\mu_1=0.6$ and $\epsilon = {4}/{3}$ when we use one layer,
and $\mu = (1.4905,  0.27594)$, $\epsilon = ( 1.09271, 1.6702)$
when we use two layers. The first two columns show the
distribution of $\mu$ and $\ep$ for $L=0,1,2$. The last column
shows the values of coefficients of $W_n[\mu, \epsilon, t]$:
$W_0^{1,0}$, $W_0^{2,0}$, $W_0^{2,1}$ indicated by $(0,0,0)$;
$W_1^{0}$, $W_1^{1,0}$, $W_1^{1,1}$ indicated by $(1,1,1)$;
$W_2^{0}$ indicated by $2$. The figures show that if $L=0$ then
none of the coefficients is zero; if $L=1$, then the coefficient
of $t^2$ is close to zero; if $L=2$ all the coefficients are close
to zero.

\begin{figure}[h!]
\begin{center}
\epsfig{figure=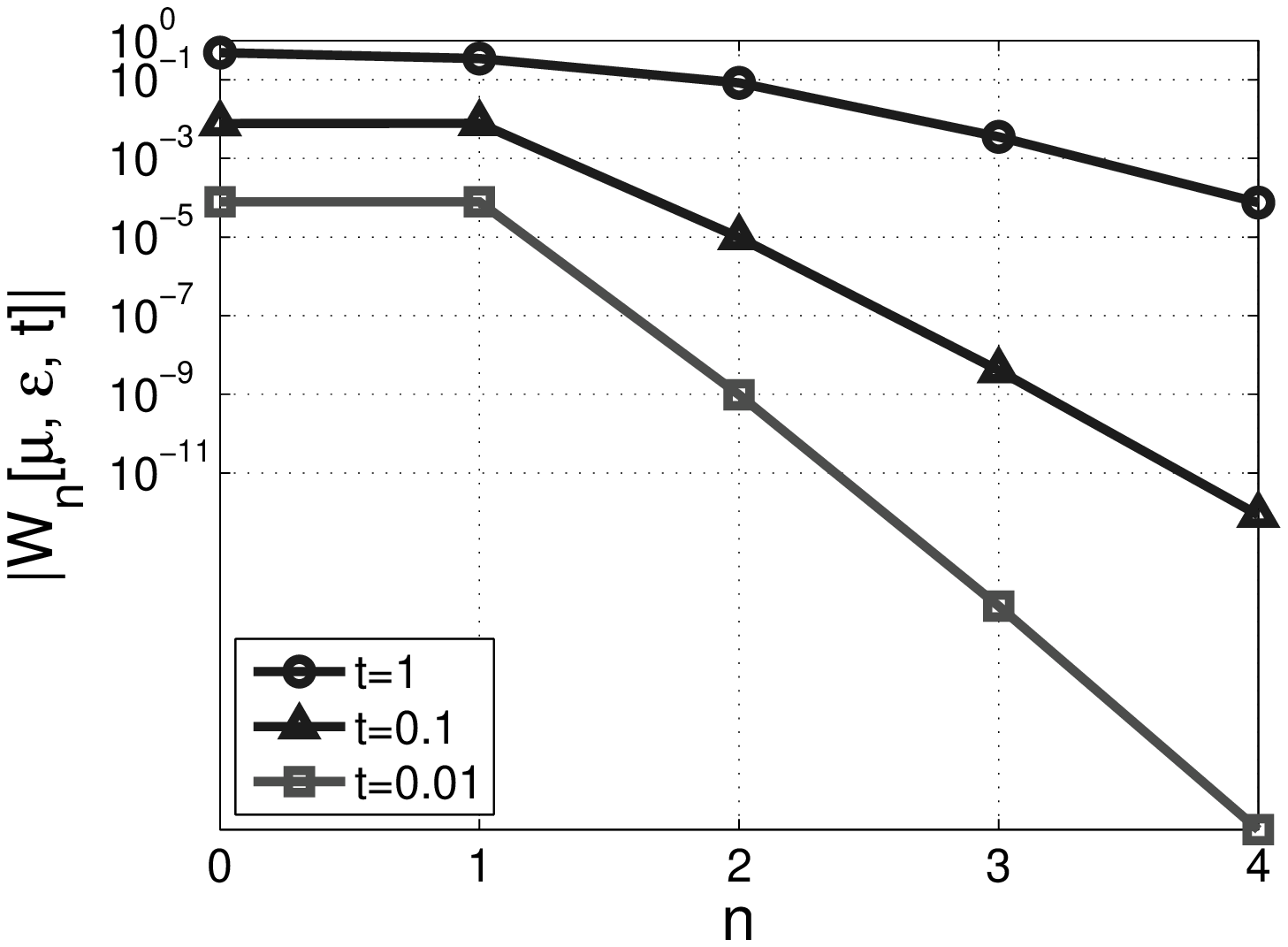,  height=5cm}\hskip 1cm
\epsfig{figure=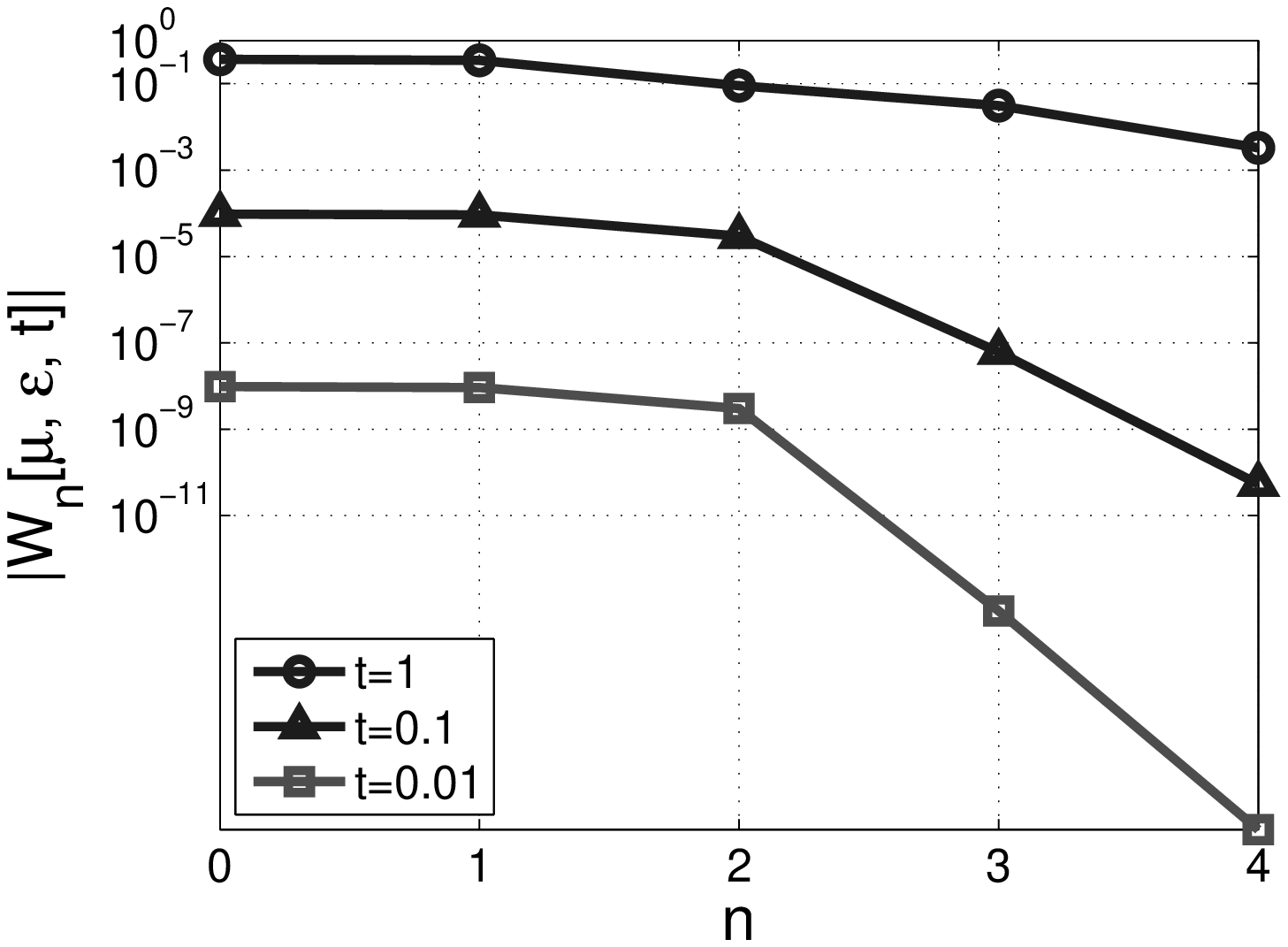,  height=5cm}
\epsfig{figure=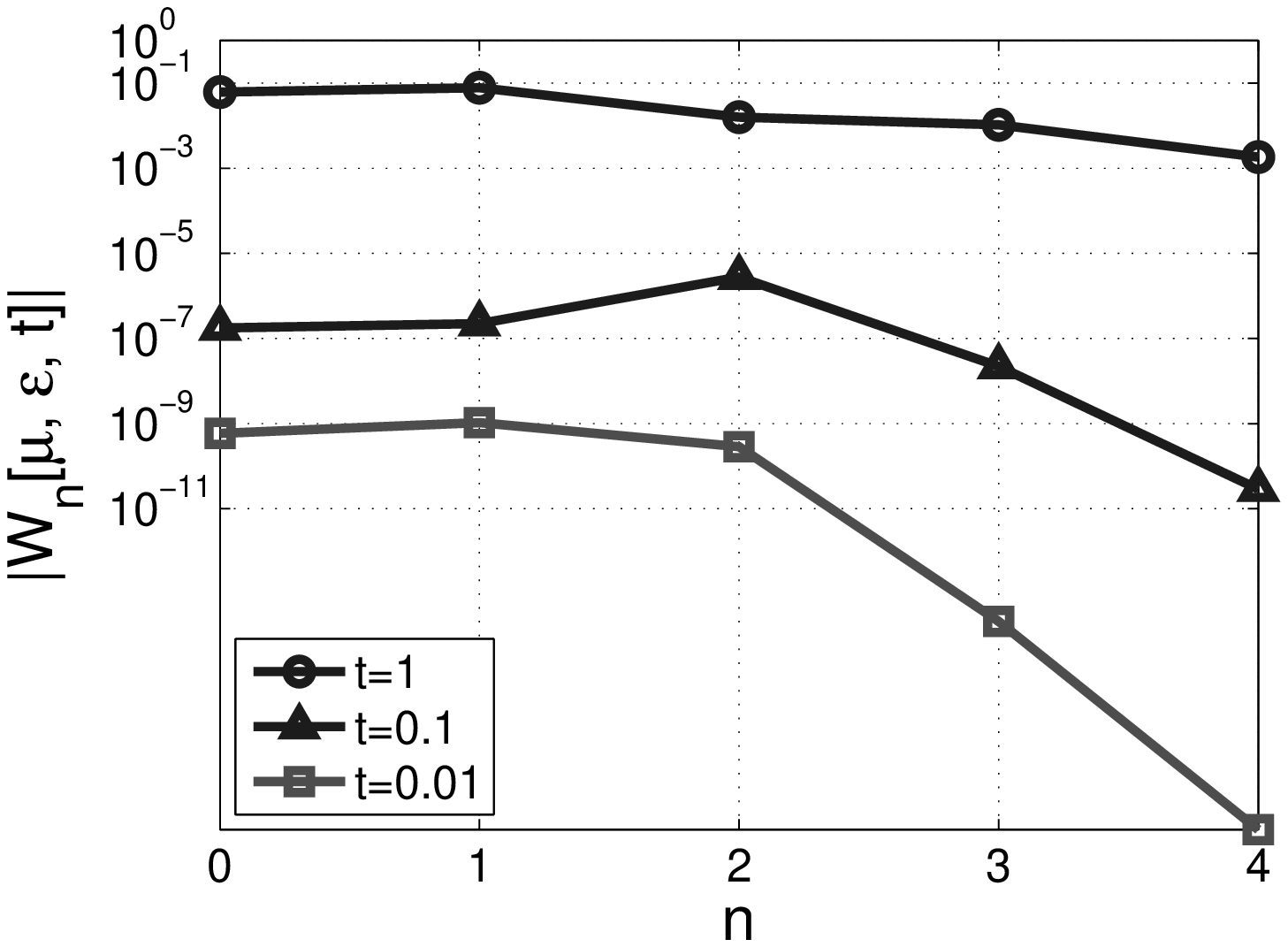,  height=5cm}\\
\end{center}
\caption{The scattering coefficient, $W_n[\mu,\ep, t]$, for
$n=0,\ldots, 4$ and $t=1, 0.1, 0.01$ using the permeability and
permittivity profiles computed in Figure \ref{general}.}
\label{general2}
\end{figure}

We then compute the scattering coefficients $W_n[\mu, \ep, t]$
using the computed $(\mu, \ep)$ for $L=0,1,2$. Figure
\ref{general2} shows the results of computations for $n=0, \ldots,
4$ and $t=1, 0.1, 0.01$. It clearly shows that $W_n[\mu, \ep, t]$
for $n \le 2$ gets smaller as the number of layers increases.

\section{Conclusion}
We have shown near-cloaking examples for the Helmholtz equation.
We have designed a cloaking device that achieves enhanced cloaking
effect. Any target placed inside the cloaking device has an
approximately zero scattering cross section. Such cloaking device
is obtained by the blow up using the transformation optics of a
multi-coated insulating inclusion. In the numerical example, to
minimize the scattering coefficients of the  multi-coated
insulating inclusion up to the second order, we set $2$ layers
 with different permittivity and
permeability properties. The technique presented in this paper can
be extended to full Maxwell's equations. This will be the subject
of a forthcoming paper.

\end{document}